\documentclass{amsart}

\usepackage{amsmath} 
\usepackage{amssymb}
\usepackage{latexsym}
\usepackage[mathscr]{eucal}
\usepackage{enumerate}
\newtheorem{theorem}{Theorem}[section]
\newtheorem{corollary}[theorem]{Corollary}

\newtheorem{lemma}[theorem]{Lemma}
\newtheorem{proposition}[theorem]{Proposition}
\newtheorem{remark}[theorem]{Remark}

 \def\Ind{\operatorname{Ind}}
 \def\diag{\operatorname{diag}}
 \def\Jord{\operatorname{Jord}}
 \def\im{\operatorname{Im}}
 \numberwithin{equation}{section}
 \def\Sym{\operatorname{Sym}}
 \def\Fr{\operatorname{Fr}}
 \def\ds{\displaystyle}
 \def\A{\bold A}
 \def\G{\bold G}
 \def\P{\bold P}
 \def\M{\bold M}
 \def\N{\bold N}
 \def\B{\bold B}
 \def\T{\bold T}
 \def\U{\bold U}
 \def\H{\bold H}
 \def\Z{\Bbb Z}
 \def\R{\mathbb R}
 \def\a{\alpha}
 \def\b{\beta}
 \def\e{\varepsilon}
 \def\k{\kappa}
 \def\l{\lambda}
 \def\o{\omega}
 \def\s{\sigma}
 \def\vp{\varphi}
 
 \def\O{\Omega}
 \def\T{\Theta}
 \def\C{\mathbb C}
 \def\Z{\mathbb Z}
 \def\fa{\mathfrak a}
 \def\fz{\mathfrak z}
 \def\sE{\mathcal E}
 
 \newcommand{\ca}{\check\alpha}

\begin{document} 
\title{$R$--groups, elliptic representations, and parameters for $GSpin$ groups}

\author{Dubravka Ban}
\author{David Goldberg}

\address{Department of Mathematics\\
Southern Illinois University\\ Carbondale, IL 62901 \\USA }

\email{dban@math.siu.edu}

\address{Department of Mathematics\\ Purdue University\\
West Lafayette, IN 47907\\ USA}

\email{goldberg@math.purdue.edu}
\maketitle
\begin{abstract}
We study parabolically induced representations  for $GSpin_m(F)$ with $F$ a $p$--adic field of characteristic zero.   The  Knapp-Stein $R$--groups  are described
and shown to be elementary two groups. We show the associated cocycle is trivial proving multiplicity one for induced representations. We classify the elliptic tempered spectrum. 
For $GSpin_{2n+1}(F)$, we describe the Arthur (Endoscopic) $R$--group attached to Langlands parameters, and show these are isomorphic to the corresponding Knapp-Stein $R$--groups.
\end{abstract}

\section*{Introduction}
We continue our study of parabolically induced representations for $p$--adic groups of classical type.  Here we turn our attention to the group $GSpin_m(F),$ as defined by Asgari \cite{asgari}. These are groups of type $B_{[m/2]}$ if $m$ is odd and type $D_{m/2}$ if $m$ is even.  A long term goal is to study the group $Spin_m(F),$ which is the simply connected split group of type $B$ or $D,$ depending on whether $m$ is odd or even, respectively.  The advantage of studying $GSpin$ groups is their Levi subgroups are nicer, making the problem more tractable.  We hope to apply the information derived here to $Spin$ groups, and we leave this to further study.

Let $F$ be a nonarchimedean field of characteristic zero, and suppose $\G$ is a connected reductive quasi-split group defined over $F.$  We denote the $F$--points, $\G(F),$ by $G$ and use this notational convention throughout this manuscript.  The admissible dual of $G$ can be be studied through the  theory of parabolically induced representations, as described in Harish-Chandra's philosophy of cusp forms \cite{hc-wt}. Moreover, the discrete, tempered, and admissible spectra are classified through parabolic induction from supercuspidal, discrete series, and tempered representations (via the Langlands Classification) \cite{silberger}.  One also wishes to divide the tempered spectrum into the elliptic classes, \cite{arthur-elliptic}, which are those which contribute to the Placherel Formula, and the non-elliptic classes.  For this purpose, we let $\sE_c(G),\sE_t(G),\sE_2(G),$ and $^\circ\sE(G)$ be the classes of irreducible  admissible, tempered, discrete series, and unitary supercuspidal representations, respectively, of $G.$  We make no distinction between a representation $\pi$ and its class $[\pi]\in\sE_c(G).$

Let $\P=\M\N$ be a parabolic subgroup of $\G,$ and suppose $\s\in\sE_2(M).$  We let $\Ind_P^G(\s)$ or $i_{G,M}(\s)$ denote the representation of $G$ obtained through normalized induction from $P,$ with $\s$ extended trivially from $M$ to $P.$  In the case of archimedean groups, Knapp and Stein developed
 the theory of standard and normalized intertwining operators (see \cite{knapp-stein}, for example).  Through a combinatorial study of the inductive properties of these normalized intertwining operators they were able to describe a finite group, $R(\s),$ whose representation theory classifies the components of $i_{G,M}(\s),$ in that there is a bijection $\rho\mapsto\pi_\rho$ from the irreducible representations $\widehat{R(\s)}$ to the inequivalent components of $i_{G,M}(\s).$  More precisely, the intertwining algebra $\mathcal C(\s)$ of $i_{G,M}(\s)$ is isomorphic to the twisted group algebra $\C[R(\s)]_{\eta},$ with $\eta$ a particular $2$--cocylce of $R(\s)$ arising from composition of intertwining operators \cite{arthur-elliptic,keys-ens}. In the archimedean case, $R(\s)$ is always abelian (in fact an elementary $2$--group), so each $\rho$ is a character and $\pi_\rho$ appears in $i_{G,M}(\s)$ with multiplicity one.  Silberger \cite{silberger-rgps, silberger-rgps2}  extended the theory of $R$-groups to $p$--adic fields.  Knapp and Zuckerman \cite{knapp-zuck-m1} showed there are cases when $R(\s)$ would be non-abelian, and hence the multiplicity of $\pi_\rho$ could be greater than one. 

If $\G=\G_n=GSpin_{2n},$ or $GSpin_{2n+1},$ then any Levi subgroup is of the form 
\begin{equation*}\label{levi-decomp}
\M\simeq GL_{n_1}\times\cdots\times GL_{n_r}\times \G_m,
\end{equation*}
with $n_1+\cdots+n_r+m=n.$  So, for any $\s\in\sE_2(M)$ we have 
\begin{equation*}\label{rep-decomp}
\s\simeq\s_1\otimes\cdots\otimes\s_r\otimes\tau,
\end{equation*}
with $\s_i\in\sE_2(GL_{n_i}(F)),$ and $\tau\in\sE_2(G_m).$  The similarity between this situation and that of the classical groups $Sp_{2n}(F)$ and $SO_n(F),$ makes it amenable to the techniques of \cite{goldberg-rgps}.  In fact we prove the $R$--groups have the same structure as these classical groups.  Thus, our first main results  can be phrased as  $R$-groups for $GSpin$ groups mirror those for split classical groups (cf. Theorems \ref{R-gp odd} and \ref{R-gp even}).  In particular, $R(\s)\simeq\Z_2^d,$ for some $0\leq d\leq r.$


Arthur \cite{arthur-elliptic} undertook the study of the elliptic spectrum, and was able to use the extension of $R(\s)$ defined by $\eta$ to characterize when components of $i_{G,M}(\s)$ have elliptic components.  Herb \cite{herb} used this characterization, along with the description of the $R$--groups in \cite{goldberg-rgps}, to determine the elliptic tempered spectrum of $Sp_{2n}(F)$ and $SO_n(F).$ Because the description of $R$--groups in our case  is similar to that of \cite{goldberg-rgps}, the techniques of \cite{herb} can be applied, and again the results are similar. To be more precise,  the cocycle $\eta$ always splits  and $i_{G,M}(\s)$ has elliptic components if and only if $d$ is as large as possible (this turns out to be $d=r$ or $r-1$ cf.  Lemma \ref{cocycle} and Theorems \ref{elliptic odd} and \ref{elliptic even}).

On the other hand, the local Langlands conjecture predicts a canonical bijection $\vp\rightarrow\Pi_\vp(G)$ between admissible homomorphisms $\vp:W_F'\rightarrow\,^LG$ and $L$--packets $\Pi_\vp(G)$ of $G.$  Here, $W_F'$ is the Weil-Deligne group, $^LG=\hat G\rtimes W_F$ is the Langlands $L$--group, with $\hat G$ the connected Langlands dual group, and $W_F$ is the Weil group.  The $L$--packets $\Pi_\vp(G)$ are finite sets which partition $\sE_c(G),$ and the members of $\Pi_\vp(G)$ are to be $L$--indistinguishable, in the sense that the Langlands $L$--functions and $\varepsilon$--factors are to be constant on $\Pi_\vp(G).$  If $\s\in\sE_2(M,)$ and $\vp:W_F'\rightarrow \,^LM$ is its Langlands parameter (i.e., $\s\in\Pi_{\vp}(M)$),  then composing with the inclusion $^LM\hookrightarrow\,^LG$ gives an $L$-packet $\Pi_{\vp}(G),$ and the elements of this $L$--packet should be all components of $i_{G,M}(\s'),$ with $\s'\in\Pi_\vp(M).$ Langlands predicted the $R$--group, $R(\s)$ should be encoded in this arithmetic information, and Arthur made this more precise in \cite{arthur-aster}. In particular, Arthur defined a finite group $R_{\vp,\s}$ attached to $\vp$ and $\sigma,$ and predicts $R(\s)\simeq R_{\vp,\s}.$  This conjecture has been confirmed in many cases \cite{ban-goldberg,B-Zh,ChLi, choiy-goldberg,goldberg-dual,  keys-ens,shelstad}.

Here we are able to prove  $R(\s)\simeq R_{\vp,\s}$ for $GSpin_{2n+1}$ in several steps.  The first is to reduce the isomorphism to the case where $\M$ is maximal, and this we do in the wider context of split groups (cf. Lemma \ref{maxcase1}).  Arthur identifies the stabilizer $W(\s)$ of $\s$ in the Weyl group with  a subgroup $W_{\vp,\s}$ of a certain Weyl group in $\hat M.$  $R(\s)$ can be realized as a quotient $W(\s)/W'$ of $W(\s),$ with $W'$ the subgroup of $W(\s)$ generated by root reflections in the zeros of the rank 1 Plancherel measures. On the other hand $R_{\vp,\s}=W_{\vp,\s}/W^\circ_{\vp,\s}$ is a quotient of $W_{\vp,\s},$ where $W_{\vp,\s}^\circ$ is the the intersection of $W_{\vp,\s}$ with another, smaller Weyl group. Thus, it is enough to show, under the isomorphism of $W(\s)$ and $W_{\vp,\s}$  that $W'$ is identified with $W_{\vp,\s}^\circ.$  Hence it is enough to show $W_{\vp,\s}^\circ$ is  generated by co-root reflections  coming from the roots for which the Plancherel measures are zero.  Shahidi  \cite{Sh2} showed, in the generic case,   the zeros of the rank 1 Plancherel measures are equivalent to poles of Langlands $L$--functions, $L(s,\s,r_i),$ (where $i=1,2$ is determined in a particular way \cite{Sh1} and $r_i$ is a representation of $^LM$ coming from its adjoint representation). The local Langlands conjecture
predicts $L(s,\s,r_i)=L(s,r_i\circ\vp),$ where the right hand side is the Artin $L$--function.  
We separate the proof of the isomorphism of  Knapp-Stein and Arthur $R$--groups into two  (maximal) cases, the Siegel parabolic subgroup, i.e
$\M\simeq GL_n\times GL_1,$ and the non-Siegel maximal parabolic subgroups, $\M\simeq GL_k\times \G_m,$ with $m\geq 2.$
The final results in these two cases can be found in Corollary \ref{corollary siegel} and 
Theorem \ref{isomorphism general case}.  For the 
latter we need  conjecture 9.4 of \cite{Sh2} (otherwise known as the Tempered $L$--packet Conjecture).   

The structure and isomorphism of Knapp-Stein and Arthur $R$--groups plays a crucial role in the transfer of automorphic forms from classical to general linear groups in \cite{arthur-book},  and among the important results therein is a proof of the Tempered $L$--packet Conjecture in the case of classical groups. We expect if the methods of \cite{arthur-book} can be extended to $GSpin$ groups, then the isomorphism of $R(\s)$ and $R_{\vp,\s}$ would play a similar role.

In Section 1 we recall the basic facts about the $GSpin$ groups.  In Section 2 we work to determine the zeros of the Plancherel measures and compute the $R$--groups for $GSpin$ groups.  In Section 3 we show the cocycle which, along with the $R$--group, determines the structure of $i_{G,M}(\s)$ is a coboundary.   We then use the results of Section 2 to classify the elliptic tempered spectra of $GSpin$ groups.   In Section 4 we prove the isomorphism of the Knapp-Stein and Arthur $R$--groups for the $GSpin_{2n+1}$ groups.

\section{Preliminaries}\label{preliminaries}
Let $F$ be a local nonarchimedean field of characteristic zero.  Let $\G=\G_n=GSpin_{2n},$ or $GSpin_{2n+1}.$  We adopt the convention that $G_0=GL_1.$We let $\H=Spin_{2n}$ or $Spin_{2n+1}.$  We recall the exact sequence
$$1\rightarrow \Z_2\rightarrow \H\rightarrow \H'\rightarrow 1,$$
where $H'=SO_{2n}$ or $SO_{2n+1},$
We have $\G$ and $\H$  are of type $D_n$ in the first case and type $B_n$ is the second case.  
Let $\hat G$ be the connected component of the Langlands $L$--group.  Then $\hat G=GSO_{2n}(\C)$ if $\G=GSpin_{2n}$ and is $GSp_{2n}$ if $\G=GSpin_{2n+1}.$  Then since $\G$ is split, $^LG=\hat G\times W_F,$ with $W_F$ the Weil group of $F.$
We fix $\B$ to be the Borel subgroup in $\G$  lying over the upper triangular Borel subgroup in $\H'.$  Let $\B=\bold T\U$ be the Levi decomposition of $\B.$  Let $\Phi=\Phi(\G,\bold T)$ be the roots of $\bold T$ in $\G,$ and let $\Delta$ be the simple roots determined by $\B.$  Then $\Delta=\{\alpha_1,\a_2,\dots,\alpha_n\},$
where  $\a_i=e_i-e_{i+1},$ for $i=1,2,\dots,n-1,$ and 
$$\alpha_n=\begin{cases}
e_{n-1}+e_n&\text{ if }\G=GSpin_{2n},\\
e_n&\text{ if }\G=GSpin_{2n+1}.
\end{cases}$$
Recall the Weyl group is $W=W(\G,\bold T)=N_\G(\bold T)/\bold T.$  Note, if $\G=GSpin_{2n+1},$ then  $W\simeq S_n\ltimes \Z_2^n,$ while if $\G=GSpin_{2n},$ we have $W\simeq S_n\ltimes\Z_2^{n-1}.$  One can compute this directly from the description in \cite{asgari-shahidi-djm}, or one can note that $W(\hat G,\hat T)$ is of this form, and use duality.  Taking this last approach,  the description of these Weyl groups given in \cite{goldberg-annalen}, which we summarize. Note
$$\hat T=\left\{\diag\left\{a_1,a_2,\cdots,a_n,\l a_n^{-1},\dots,\l a_2^{-1},\l a_1^{-1}\right\}\big | a_i,\l\in\C^\times\right\}$$
in either case.  We may denote an element of $\hat T$ by $t(a_1,a_2,\dots,a_n,\l).$  If $s\in S_n,$ then we also denote by $\hat s$  a representative of the element of $W(\hat G,\hat T)$ such that $\hat s t(a_1,a_2,\dots,a_n,\l)\hat s^{-1}=t(a_{s(1)},a_{s(2)},\dots,a_{s(n)},\l).$
If $\G=GSpin_{2n+1},$ then denote by $\hat c_i$ a representative of the element of $W(\hat G,\hat T)$ such that $\hat c_i t(a_1,\dots,a_i,\dots, a_n,\l)\hat c_i^{-1}=t(a_1,\dots,\l a_i^{-1},\dots,a_n,\l).$  Then $W(\hat G,\hat T)$ is generated by $\{\hat s|s\in S_n\}$ and $\{\hat c_i|1\leq i\leq n\}.$
If $\G=GSpin_{2n},$ then $W(\hat G,\hat T)$ is generated by $\{\hat s|s\in S_n\}$ and $\{c_ic_j|1\leq i,j\leq n\}.$
We have the pairing of roots $\Phi(\G,\bold T)$ and coroots $\Phi(\hat G,\hat T)$ which we denote by $\a\mapsto\check\a,$  and denote by $w$ the element of $W(\G,\bold T)$ corresponding to $\hat w$ by this pairing.

Let $\P=\M\N\supset \B$ be a standard parabolic subgroup of $\G.$  
Then, for some  $\theta\subset\Delta$ we have $\P=\P_\theta=\M_\theta\N_\theta.$
Then  there is a partition $n=n_1+n_2+\dots+n_r+m,$ so that 
 $\theta=\Delta\setminus\{\a_{n_1},\a_{n_1+n_2,},\dots,\a_{n_1+n_2+\dots +n_r},\a_n\},$ if $m=0,$ and
 $\theta=\Delta\setminus\{\a_{n_1},\a_{n_1+n_2,},\dots,\a_{n_1+n_2+\dots +n_r}\},$ if $m>0.$
 Then 
 \begin{equation}\label{levi} 
 \M\simeq GL_{n_1}\times GL_{n_2}\times \dots\times GL_{n_r}\times \G_m.
 \end{equation}
 Let $\A$ be the split component of $\P,$ and let $\Phi(\P,\A)$ be the reduced roots of $\A$ in $\P.$
For $i=1,2,\dots,r,$ we let $k_i=n_1+\dots+n_i.$  
 Then, for $1\leq i < j\leq r,$  set $\a_{ij}=e_{k_i}-e_{{k_{j-1}}+1},$ and $\b_{ij}=e_{k_i}+e_{{k_{j-1}}+1},$ and 
$$\gamma_i=\begin{cases} e_{k_i}+e_n&\text{ if } \G =GSpin_{2n};\\ e_{k_i}& \text{ if }\G=GSpin_{2n+1}.\end{cases}$$

We describe the relative Weyl group $W_\M=N_\G(\A_M)/Z_\G(\A_M)=N_{\G}(\A_{\M})/\M.$  Suppose $\M$ is as above.  As in the case of other groups of classical type, $W_M\subset S_r\ltimes \Z_2^r.$  If $\G$ is of type $B_n,$ then $W_\M\simeq S\ltimes\Z_2^r,$ for some subgroup $S$ of $S_r.$  In fact $\ds{S=\langle(ij)| i<j, n_i=n_j\rangle.}$ 
More precisely, let $k_0=0,$ and  for $i=1,2,\dots,r-1,$  let $k_i$ be as above.  If $n_i=n_j,$ let $[ij]\in W(\G,\bold T)$ be the element 
$\ds{\prod_{k=1}^{n_i}(k_{i-1}+k\, k_{j-1}+k).}$  Then $[ij]\mapsto (ij)$ gives an isomorphism of $W_M\cap S_n$ to $S.$  We generally denote these elements by the more standard $(ij).$ For $1\leq i\leq r,$ we let 
$\ds{C_i=\prod_{k=1}^{n_i}c_{k_{i-1}+k}.}$  We call $C_i$ a {\bf block sign change}, and $\langle C_i|i=1,\dots,r\rangle \simeq\Z_2^r$ is the sign change subgroup of $W_\M.$
The action of $S$ on $\M$ is given by 
$$(ij):(g_1,\dots,g_r,h)=(g_1,\dots,g_{i-1},g_j,g_{i+1},\dots,g_{j-1},g_i,\dots,g_r,h).$$  
Also, from the action of $C_i$ on the root datum of $\G$ (see \cite{asgari}) we have $C_i\cdot (g_1,\dots,g_i,\dots,g_r,h)=(g_1,\dots,\!^tg_i^{-1},\dots,g_r,e_0^*(\det g_i)h).$
If $\G$ is of type $D_n,$ then $W_\M\simeq S\ltimes \mathcal C,$ where $S$ is as above for type $B_n,$ and $\mathcal C\subset\Z_2^r.$  If $m=0,$ then we  have $\mathcal C=\mathcal C_1\times\mathcal C_2,$ where $\mathcal C_1=\langle C_i|n_i\text{ is even}\rangle,$ and $\mathcal C_2=\langle C_iC_j|n_i,n_j\text{ are odd}\rangle.$  If $m>0,$ then
$\mathcal C\simeq\Z_2^r,$  and $$\mathcal C=\langle C_i|n_i\text{ is even}\rangle\times \langle C_ic_n|n_i\text{ is odd}\rangle.$$
We note that $S$ and each $C_i$ acts as in the case of type $B_n,$   (and of course $C_iC_j$ acts as the product in type $D_n$). In the case of $m>0$ and $n_i$ odd, we have
$C_ic_n\cdot(g_1,\dots, g_i,\dots, g_r,h)=(g_1,\dots,^tg_i^{-1},\dots, g_r,(\det g_i)(c_n\cdot h)),$
where $c_n$ is given by the outer automorphism on the Dynkin diagram of $\G_m.$

\section{R-groups for GSpin}
We continue with the notation of the previous section.  Let $\M$ be a Levi subgroup of $\G=\G_n$ and assume $\M$ is of the form \eqref{levi}. Let $\s\in\mathcal E_2(M).$ Then
$\s\simeq\s_1\otimes\s_2\dots\otimes\s_r\otimes\tau,$
where
$\s_i\in\mathcal E_2(GL_{n_i}(F)),$ and $\tau\in\mathcal E_2(G_m).$
For $\a\in\Phi(\P,\A),$ we set $\A_\a=(\A\cap\ker\a)^\circ,$ and $\M_\a=Z_\G(\A_\a).$   Then $^*\P_\a=\P\cap\M_\a=\M\N_\a,$
where $\N_\a=\N\cap \M_\a$ is a maximal parabolic subgroup of $\M_\a$ with Levi component $\M.$   We let $W_\a=W(\M_\a,\A).$  If $W_\a\neq\{1\},$ we let $w_\a$ be the unique nontrivial element of $W_\a.$
We recall the Plancherel measure, $\mu_\a(\s)$ is determined by the standard intertwining operator attached to $\Ind_{^*\!P_\a}^{M_\a}(\s),$ and in particular, $\mu_\a(\s)=0$ if and only if $w_\a\s\simeq\s$ and $\Ind_{^*\!P_\a}^{M_\a}(\s)$ is irreducible.

We note if $\a=\a_{ij},$ then 
\begin{equation}\label{2.1}
M_\a\simeq\prod_{k\neq i,j} GL_{n_k}\times GL_{n_i+n_j}\times\G_m,
\end{equation}
 and 
\[
W_\a=\begin{cases} 1&\text{ if }n_i\neq n_j;\\\{1,(ij)\}&\text{ if } n_i=n_j.\end{cases}
\]
If $\a=\b_{ij},$ then 
$\M_\a\simeq\M_{\a_{ij}}$ is again given by \eqref{2.1}, and 
$$W_\a=\begin{cases} 1&\text{ if }n_i\neq n_j;\\\{1,(ij)C_iC_j\}&\text{ if } n_i=n_j.\end{cases}$$
Finally, for $\a=\gamma_i,$ we have 
$$\M_\a\simeq\prod_{k\neq i}GL_{n_k}\times\G_{n_i+m}.$$
If $\G$ is of type $B_n,$ or $n_i$ is even, then $W_\a=\{1,C_i\}.$  If $\G$ is of type $D_n,$ and $n_i$ is odd,
then $$W_\a=\begin{cases} C_ic_n&\text{ if } m>0;\\\text{1}&\text{ if } m=0.\end{cases}$$

 We note, for $\G$ of type $B_n,$
$$w_\a\s\simeq\begin{cases} \s_1\otimes\cdots\otimes\s_{i-1}\otimes\s_j\otimes\s_{i+1}\cdots\otimes\s_{j-1}\otimes\s_i\otimes\cdots\otimes\s_r\otimes\tau;\\\s_1\otimes\cdots\otimes\s_{i-1}\otimes(\tilde\s_j\otimes\o_\tau)\otimes\s_{i+1}\cdots\s_{j-1}\otimes (\tilde\s_i\otimes\o_\tau)\otimes\s_{j+1}\otimes\cdots\otimes\s_r\otimes\tau;\\
\s_1\otimes\cdots\otimes\s_{i-1}\otimes(\tilde\s_i\otimes\o_\tau)\otimes\cdots\otimes\s_r\otimes\tau,\end{cases}$$
if $\s=\a_{ij},\b_{ij},$ or $\gamma_i,$ respectively.  Here $\o_\tau$ is the central character of $\tau$ restricted to the identity component of the center of $G_m.$ 
For  type $D_n,$  the result is as above,  except in the case where $\a=\gamma_i,$ $n_i$ is odd and $m>0,$ in which case
$$ w_\a \s \simeq\s_1\otimes\cdots\otimes(\tilde\s_i\otimes\o_\tau)\otimes\cdots\otimes\s_r\otimes(c_n\cdot\tau).$$

\begin{lemma} For $1\leq i <j\leq r-1$ we have $\Ind_{^*\!P_{\a_{ij}}}^{M_{\a_{ij}}}(\s)$ is irreducible.  Similarly $\Ind_{^*\!P_{\b_{ij}}}^{M_{\b_{ij}}}(\s)$ is irreducible.
\end{lemma}
\begin{proof} Let $\a=\a_{ij}.$
In this case $M_\a$ is given by \eqref{2.1}.
Let $\bold Q_{ij}$ be the standard $GL_{n_i}\times GL_{n_{j}}$--parabolic subgroup of $GL_{n_i+n_{j}}.$  Then,
$$\Ind_{^*\!P_\a}^{M_\a}(\s)\simeq\left(\bigotimes_{\ell\neq i,j+1}\s_\ell\right)\otimes\left( \Ind_{Q_{ij}}^{GL_{n_i+n_{j}}(F)}(\s_i\otimes\s_{j})\right)\otimes\tau,$$ and the result now follows from Olsanskii or Bernstein and Zelevinski\cite{bz,olsanskii}.

If $\a=\b_{ij},$ then we again have $M_\a$ is given by \eqref{2.1},
and in this case
$$\Ind_{^*\!P_\a}^{M_\a}(\s)\simeq\left(\bigotimes_{\ell\neq i,j}\s_\ell\right)\otimes\left( \Ind_{Q_{ij}}^{GL_{n_i+n_{j}}(F)}(\s_i\otimes(\tilde\s_j\otimes\o_\tau))\right)\otimes\tau.$$  Thus, the result again follows from \cite{bz,olsanskii}.
\end{proof}

From this we derive the following result.

\begin{corollary}\label{planch=0}
If $\a=\a_{ij},$  then $\mu_\a(\s)=0$ if and only if $n_i=n_j$ and $\s_i\simeq\s_j.$
If $\a=\b_{ij},$  then $\mu_\a(\s)=0$ if and only if $n_i=n_{j}$ and $\s_i\simeq\tilde\s_{j}\otimes\o_\tau.$
\end{corollary}

\begin{lemma} Let $\sigma=\sigma_1\otimes\s_2\otimes\cdots\otimes\s_r\otimes\tau\in\mathcal E_2(M),$ and let $R=R(\s).$  Suppose $w\in R$ and $w=sc,$ with $s\in S_r$ and $c\in \Z_2^r.$  Then $s=1.$
\end{lemma}

\begin{proof} This is a {\bf Keys argument} as defined in \cite{goldberg-rgps} and introduced in \cite{keys}. Since the sign changes act independently on the disjoint cycles of $s,$
we may suppose, without loss of generality, that $s=(12\cdots j)$.  Furthermore, if $\G$ is of type $B_n,$ then
up to conjugation by sign changes we may assume $c=C_jc',$ or $c=c',$with $c'$ not changing signs among $1,2,\dots,j.$   If $\G$ is of type $D_n,$ then we may assume $c$ is either of the same form, or of the form $C_{j-1}C_jc',$  with $c'$ changing no signs among $1,2,\dots,j.$
 If $c$ changes no (block) signs among $1,2,\dots,j,$ then we note that $\s_1\simeq\s_2\simeq\cdots\simeq\s_j.$  So, in particular $\a_{1j}\in\Delta'$ and $w(\a_{1j})=-\a_{12}<0,$ so $w\not\in R(\s).$  If $c=C_jc',$ then $\s_1\simeq\s_2\simeq\cdots\simeq\s_{j-1}\simeq\s_j\simeq(\tilde\s_1\otimes\o_\tau)\,$ and thus $\b_{1j}\in\Delta'.$  However, $w\b_{1j}=-\a_{12}<0,$ so $w\not\in R.$  
 Finally, if $c=C_iC_j c',$ then $w\s\simeq\s$ implies $\s_1\simeq\s_2\simeq\cdots\simeq\s_{j-1}\simeq(\tilde\s_j\otimes\o_\tau), $ and therefore, again, $\b_{1j}\in\Delta',$ with $w\b_{1j}=-\a_{12}<0, $ showing $w\not\in R.$
  \end{proof}

\begin{corollary}\label{r-is-2-gp} For $G=GSpin_{2n}$ or $GSpin_{2n+1},$ we have $R\subset \Z_2^r.$ 
\end{corollary}

We let $W(\s)=\{w\in W(\G,\A_\M)|w\s\simeq\s\}.$  If $\G$ is of type $B_n,$ and  $W(\s)\neq 1,$ then one of the following holds:

\begin{align}
 \s_i\simeq&\s_j, \text{ for some }i\neq j;\label{2.2}\\
 \s_i\simeq&\tilde\s_j \otimes\o_\tau\text{ for some }i\neq j; and\label{2.3} \\
\s_i\simeq&\tilde\s_i\otimes\o_\tau.\label{2.4}
 \end{align}
 
Note that  \eqref{2.2} holds if $(ij)\in W(\s),$  \eqref{2.3} holds if $(ij)C_iC_j\in W(\s),$ while  \eqref{2.4} holds if $C_i\in W(\s).$ Also notice if $w=(ij)C_i\in W(\s),$ then $w^2=C_iC_j\in W(\s),$ so this case is covered by \eqref{2.4}.  
For $w\in W(\G,\A),$ we let $R(w)=\{\a\in\Phi(\P,\A)|w\a<0\}.$

For $B\subset\{1,2,\dots,r\},$ we let $\ds{C_B=\prod_{i\in B}C_i.}$
If $C_B\in R(\s),$ then $R(C_B)\cap\Delta'=\emptyset.$
Note that 
 $$R(C_B)=\left\{\a_{ij},\b_{ij}\Big | i\in B, i<j\right\}\cup\left\{\gamma_i|i\in B\right\}.$$
We let $\bold Q_i$ be the standard $GL_{n_i}\times \G_m$ parabolic subgroup of $\G_{n_i+m}.$

 \begin{theorem}\label{R-gp odd} Let $\G=GSpin_{2n+1}$
and $\M\simeq GL_{n_1}\times\cdots\times GL_{n_r}\times \G_m,$ with $\ds{m+\sum_i n_i=n.}$
Let $\s\simeq\s_1\otimes\cdots\otimes\s_r\otimes\tau\in\mathcal E_2(M),$ with each $\s_i\in \mathcal E_2(GL_{n_i}(F))$ and $\tau\in\mathcal E_2(G_m).$  Let $d$ be the number of nonequivalent classes among $\{\s_1,\dots,\s_r\}$ for which $\Ind_{Q_i}^{G_{n_i+m}}(\s_i\otimes\tau)$ is reducible.  Then $R(\s)\simeq\Z_2^d.$  More precisely,  let 
$$\O(\s)=\left\{i\big |\Ind_{Q_i}^{G_{n_i+m}}(\s_i\otimes\tau)\text{ is reducible, and }\s_j\not\simeq\s_i\text{ for all }j>i\right\}.$$
Then $R(\s)=\left\langle C_i\right\rangle_{i\in\O(\s)}.$
 \end{theorem}
\begin{remark}
By Bruhat Theory we know if  $\Ind_{Q_i}^{G_{n_i+m}}(\s_i\otimes\tau)$  is reducible implies $C_i\in W(\s),$ so $\s_i\simeq\tilde\s_i\otimes\o_\tau.$
\end{remark}

\begin{proof}
From Corollary \ref{r-is-2-gp} we know $R\subset\left\langle C_i\right\rangle_{i=1}^r\simeq\Z_2^r.$
Suppose $B\subset\{1,2,\dots,r\},$ with $C_B\in R(\s).$  Then $C_B\in W(\s),$ so $\s_i\simeq\tilde\s_i\otimes\o_\tau,$ for all $i\in B.$  Thus, for each $i\in B,$ we have $C_i\in W(\s).$
Since $R(C_i)\subset R(C_B),$ and $R(C_B)\cap\Delta'=\emptyset,$ we have $R(C_i)\cap\Delta'=\emptyset.$  So $C_i\in R(\s).$  Therefore, for some subset, $\O,$ of $\{1,2,\dots,r\}$ we have $R(\s)=\langle C_i|i\in \O\rangle.$  Now suppose $C_i\in R(\s).$  For each $j>i,$ we have $\a_{ij}\in R(C_i),$ and thus $\a_{ij}\not\in\Delta'.$ By Corollary \ref{planch=0} this implies $\s_j\not\simeq\s_i,$ for all $j>i.$  Also, for each $j>i,$ we have $\b_{ij}\in R(C_i),$ so $\s_j\not\simeq\tilde\s_i\otimes\o_\tau.$  However, since $\s_i\simeq\tilde\s_i\otimes\o_\tau,$ we see $\b_{ij}\not\in\Delta'$ imposes no further condition.  Finally, since $\gamma_i\in R(C_i),$ we must have $\gamma_i\not\in\Delta'.$  We note 
$$\M_{\gamma_i}\simeq\prod_{j\neq i}GL_{n_j}\times \G_{n_i+m},$$
and 
$$\Ind_{^*\!P_{\gamma_i}}^{M_{\gamma_i}}\s\simeq\bigotimes_{j\neq i}\s_j\otimes\Ind_{Q_i}^{G_{n_i+m}}(\s_i\otimes\tau).$$
Since $C_i\in W(\s)\cap W_{\gamma_i},$ we have $\gamma_i\not\in\Delta'$ if and only if $\Ind_{Q_i}^{G_{n_i+m}}(\s\otimes\tau)$ is reducible. Thus, $i\in\O(\s),$ so $\O\subset\O(\s).$  Conversely, if $i\in\O(\s),$ then $C_i\s\simeq\s,$ and $R(C_i)\cap\Delta'=\emptyset,$ so $C_i\in \O.$  Thus $\O=\O(\s),$ and $R(\s)$ has the form we claim.
\end{proof}

Now suppose $\G$ is of type $D_n.$ Let $\M\simeq GL_{n_1}\times\cdots\times GL_{n_r}\times\G_m.$  We may assume $n_i$ is even for $i=1,2,\dots,t,$ and
$n_i$ is odd for $i=t+1,\dots,r.$   If $m=0,$ then 
$$\mathcal C\simeq \begin{cases} \Z_2^{r-1}&\text{if }t<r;\\\Z_2^r&\text{otherwise.}\end{cases}$$
If $m>0,$ then
$\mathcal C\simeq\Z_2^r,$ as described above.  
If $m=0$ or $c_n\tau\not\simeq\tau,$ then the following describes the conditions under which $W(\s)\neq\{1\}:$
\begin{align}
\s_i&\simeq\s_j \text{ for some }i\neq j;\label{2.5}\\
\s_i&\simeq\tilde\s_j\otimes\o_\tau \text{ for some }i\neq j;\label{2.6}\\
\s_i&\simeq\tilde\s_i\otimes\o_\tau\text{ for some }i\text{ with }n_i \text{ even};\label{2.7}\\
\s_i&\simeq\tilde\s_i\otimes\o_\tau \text{ and }\s_j\simeq\tilde\s_j\otimes\o_\tau\text { for some }i\neq j\text { with } n_i,n_j \text{ odd}.\label{2.8}
\end{align}
We have \eqref{2.5} holding if and only if $(ij)\in W(\s),$ \eqref{2.6} holds if and only if $(ij)C_iC_j\in W(\s),$ while \eqref{2.7} and \eqref{2.8} are the conditions for either $C_i$ (for $n_i$ even) or $C_iC_j$ to be in $W(\s).$  If $m>0$ and $c_n\tau\simeq\tau,$ then \eqref{2.5},\eqref{2.6}, and \eqref{2.7} are the conditions, with the restriction on parity removed from \eqref{2.7}.

\begin{theorem}\label{R-gp even} Let $\G=Gspin_{2n},$ and $\M\simeq GL_{n_1}\times\cdots\times GL_{n_r}\times \G_m,$ with $\ds{m+\sum_in_i=n.}$  Let $\s\in\mathcal E_2(M),$ with each $\s_i\in\mathcal E_2(GL_{n_i}(F)),$ and $\tau\in\mathcal E_2(G_m).$
\begin{enumerate}[(i)]
\item If $m=0$  or $c_n\tau\not\simeq\tau,$ then we let 
\begin{gather*}
\O_1(\s)=\{1 \leq i\leq r| n_i \text{ is even, }\Ind_{Q_i}^{G_{n_i+m}}(\s_i\otimes\tau)\text{ is reducible, and }\s_j\not\simeq\s_i\text{ for all }i>j\},\\
\intertext{and}
\O_2(\s)=\{1 \leq i\leq r| n_i \text{ is odd, }\s_i\simeq\tilde\s_i\otimes\o_\tau,\text { and }\s_j\not\simeq\s_i,\text{ for all } j>i\}.
\end{gather*}
We set $d_i=|\O_i(\s)|,$ for $i=1,2.$ Then $R(\s)\simeq\Z_2^{d_1+d_2-1},$ unless $d_2=0,$ in which case $R(\s)\simeq\Z_2^{d_1}.$  More precisely,
$$R(\s)=\left\langle C_i\big|i\in\O_1(\s)\right\rangle\times\left\langle C_iC_j\big|i,j\in\O_2(\s)\right\rangle.$$
\item If $m>0$ and $c_n\tau\simeq\tau,$ we let 
$$\O(\s)=\{1\leq i\leq r| \Ind_{Q_i}^{G_{n_i+m}}(\s_i\otimes\tau)\text{is reducible, and }\s_j\not\simeq\s_i\text{ for all }j>i\}.$$
Let $d=|\O(\s)|.$  Then $R(\s)\simeq\Z_2^d,$ and in particular,
$$R(\s)=\langle C_i|i\in\O(\s)\text{ and } n_i \text{ is even}\rangle\times\langle C_ic_n|i\in\O(\s)\text{ and } n_i\text{ is odd}\rangle.$$
\end{enumerate}
\end{theorem}

\begin{proof}
We assume $n_i$ is even for $i=1,2,\dots,t,$ and $n_i$ is odd for $i=t+1,\dots,r.$
 Suppose $m=0.$ Then 
$W_{\M}= S\ltimes\mathcal C,$
where 
$$\mathcal C=\langle C_i|1\leq i\leq t\rangle \times\langle C_iC_j| t+1\leq i,j\leq r\rangle.$$
 By Corollary \ref{r-is-2-gp}, $R(\s)\subset\mathcal C.$
Suppose $B\subset\{1,2,\dots r\}.$  Then we let $B_1=B\cap\{1,2,\dots,t\},$ and $B_2=B\setminus B_1.$
Suppose $\ds{C_B=\prod_{i\in B}C_i\in R(\s)}.$  Then $\s_i\simeq\tilde\s_i\otimes\o_\tau,$ for each $i\in B.$
Thus, $C_i\in W(\s),$ for each $i\in B_1,$ and $C_iC_j\in W(\s),$ for each $i,j\in B_2.$
As in the case of type $B_n,$ we have, for each $i\in B, \,R(C_i)\subset R(C_B),$ 
and thus $C_i\in R(\s),$ for each $i\in B_1,$ and $C_iC_j\in R(\s)$ for each $i,j\in B_2.$
Thus, there is some $\O\subset\{1,\dots,r\},$ for which
$$R(\s)=\langle C_i|\, i\in \O_1\rangle\times\langle C_iC_j|\, i,j\in \O_2\rangle.$$
For $1\leq i\leq t,$ we have
\[
R(C_i)=\{\gamma_i\}\cup\{\a_{ij},\b_{ij}\}_{j>i}.
\]
We have $R(C_i)\cap\Delta'=\emptyset,$ so by Corollary \ref{planch=0} $\s_j\not\simeq\s_i$ for all $j>i,$ as in the case of type $B_n.$
Further note, since $C_i\in W(\s),$ we have $\gamma_i\in\Delta'$ if and only if $\Ind_{^*\!P_{\gamma_i}}^{M_{\gamma_i}}\s$ is irreducible.
Since
$$\Ind_{^*\!P_{\gamma_i}}^{M_{\gamma_i}}\s\simeq\left(\bigotimes_{j\neq i}\s_i\right) \otimes\Ind_{Q_i}^{G_{n_i+m}}(\s_i\otimes\tau),$$
we see  $C_i\in R(\s)$ implies $\Ind_{Q_i}^{G_{n_i+m}}(\s_i\otimes\tau)$ is reducible.    Thus, $i\in\O_1(\s).$  Therefore, we have $\O_1\subset\O_1(\s).$  However, the reverse inclusion is now obvious.

Now suppose  $i,j\geq t+1,$ and
$C_iC_j\in R(\s).$ Then we have noted $\s_i\simeq\tilde\s_i\otimes\o_\tau,$ and $\s_j\simeq\tilde\s_j\otimes\o_\tau.$
Note further,
$$R(C_iC_j)=\{\gamma_i,\gamma_j\}\cup\{\a_{ik},\b_{ik}\}_{k>i}\cup\{\a_{j\ell},\b_{j\ell}\}_{\ell>j}.$$
As above, this now implies $\s_i\not\simeq\s_k,$ for $k>i,$ and $\s_j\not\simeq\s_\ell,$ for $\ell>j.$
Thus, we see $i,j\in\O_2(\s),$ so $\O_2\subset\O_2(\s).$
For the opposite inclusion we note, $W_{\M_{\gamma_i}}=\{1\}=W_{\M_{\gamma_j}},$ and hence $\gamma_i,\gamma_j\not\in\Delta'.$  Thus,
if $i,j\in\O_2(\s),$ then $C_iC_j\in R(\s).$  Therefore, $R(\s)$ has the form we claim.  

If $m>0$ and $c_n\tau\not\simeq\tau,$ then the argument above is essentially valid with the following adjustments. We note
$W_\M=S\ltimes\mathcal C,$
with
\begin{equation}\label{C for last case}
\mathcal C=\langle C_i|1\leq i\leq t\rangle\times\langle C_ic_n| i>t\rangle,
\end{equation}
and since $c_n\tau\not\simeq\tau,$ we have $C_ic_n\not\in W(\s),$ for $i>t.$
Also, we note for $i>t,$ $W_{\M_{\gamma_i}}=\{1,C_ic_n\},$ so $W_{\M_{\gamma_i}}\cap W(\s)=\{1\},$ and again we must have $\gamma_i\not\in\Delta'.$

(ii) Now suppose $m>0$ and $c_n\tau\simeq\tau.$  We still have 
$W_\M=S\ltimes\mathcal C,$ with $\mathcal C$ given by \eqref{C for last case}.
For $i=1,2,\dots, r,$ we let 
$$\bar C_i=\begin{cases} C_i&\text{if }i\leq t;\\C_ic_n&\text{if } i>t.\end{cases}.$$
If $B\subset \{1,2,\dots,r\},$ and
$\bar C_B=\prod_{i\in B}\bar C_i\in R(\s),$
then $\s_i\simeq\tilde\s_i\otimes\o_\tau,$ for each $i\in B.$ So $\bar C_i\in W(\s),$ for each $i
\in B.$
Further,
$$R(\bar C_B)=\bigcup_{i\in B} R(\bar C_i),$$
so  $\bar C_i\in R(\s)$ for each $i\in B.$
Thus, there is some $\O\subset\{1,2,\dots,r\}$ such that  $R(\s)=\left\langle \bar C_i\mid i\in\O\right\rangle.$
Since
$$R(\bar C_i)=\{\a_{ij},\b_{ij}\}_{j>i}\cup\{\gamma_i\},$$
and, given $\bar C_i\in W(\s),$  we have $\a_{ij},\b_{ij}\in\Delta'$ if and only if $\s_i\simeq\s_j.$
Further, as above, $\gamma_i\in\Delta'$ if and only if $\bar C_i\in W(\s),$ and $\Ind_{Q_i}^{G_{n_i+m}}(\s\otimes\tau)$ is irreducible.
Thus,
$$\O=\{i|\Ind_{Q_i}^{G_{n_i+m}}(\s_i\otimes\tau)\text{ is reducible, and }\s_j\not\simeq\s_i,\text{ for all } j>i\}=\O(\s),$$
as claimed.
\end{proof}

\section{Elliptic Representations for $GSpin$ groups}
We now consider the question of which tempered representations of $G=GSpin_n(F)$ are elliptic.  We  can adapt the arguments of \cite{herb} 
to our current situation.   We let $G_e$ be the  set of regular elliptic elements of $G.$  If $\pi$ is an irreducible representation of $G,$ then we denote by $\T_\pi$ its character.  By Harish-Chandra \cite{HC} we know $\T_\pi$ is given by a locally integrable function, also denoted $\T_\pi,$ on the regular set.  We let $\T^e_\pi$ be the restriction of $\T_\pi$ to $G_e.$  Then $\pi\in\mathcal E_t(G)$ is elliptic if $\T_\pi^e\neq0.$

We begin by showing the $2$--cocyle arising from constructing self intertwining operators in $\mathcal C(\sigma)$ is a coboundary.
Let $\G_n=GSpin_{2n}$ or $GSpin_{2n+1}.$  Suppose $\M\simeq GL_{n_1}\times\cdots\times GL_{n_r}\times\G_m$ is a proper Levi subgroup of $\G.$  Let $\s\simeq\s_1\otimes\s_2\otimes\cdots\otimes\s_r\otimes\tau$ be an irreducible discrete series of $M.$  Let $V$ be the space of the representation $\s.$ For each $w\in R(\s),$ we choose an intertwining operator $T_w:V\rightarrow V$ so that $T_w\circ w\s=\s\circ T_w.$

\begin{lemma}\label{cocycle}  We can choose the operators $T_w$ so that $T_{w_1w_2}=T_{w_1}T_{w_2}.$
\end{lemma}

\begin{proof} For each $i,$ we let $V_i$ be the space of the representation $\s_i.$  So $V=V_1\otimes\cdots V_r\otimes V_\tau.$ Denote by $\s_i^*$ the representation on $V_i$ given by $\s_i^*(g)=\s_i(\,^tg^{-1}).$  By the work of Gelfand and Kazhdan \cite{GK} we know $\s_i^*\simeq\tilde\s_i.$
Let $\mathcal B(\s)=\{i|\s_i\simeq\tilde\s_i\otimes\o_\tau\}.$
For each $i\in\mathcal B(\s),$ we choose an intertwining operator $T_i:V_i\rightarrow V_i,$ with $T_i(\s_i^*\otimes\o_\tau)=\s_i T_i.$ 
We note $T_i^2$ is a scalar on $V_i,$ and so we can choose $T_i$ so that $T_i^2=1.$
Extend this to an operator on $V,$ by setting $T_i^V$ to be trivial on each factor, except for $V_i,$ where it is $T_i.$  
Now $T_i^V\circ C_i\s=\s T_i^V,$ and $(T_i^V)^2=\operatorname{Id}.$  Also note, for $i\neq j,$ we have $T_i^VT_j^V=T_j^VT_i^V.$
If $\G$ is of type $D_n,$ and $c_n\tau\simeq\tau,$ we choose $T_\tau$ intertwining $\tau$ and $c_n\tau,$ again with $T_\tau^2=\operatorname{Id}.$
Extend $T_\tau$ to $V$ by setting $T_\tau^V$ to be trivial on each $V_i$ and to be $T_\tau$ on $V_\tau.$
Suppose $B\subset\mathcal B(\s),$ and that 
$$w=C_B=\prod_{i\in B} C_i\in R(\s).$$
Then, we set 
$$T_w=\prod_{i\in B} T_i^V.$$
In the case where $\G$ is of type $D_n$ and $c_n\tau\simeq\tau,$ we may have 
$$w=\bar C_B=\left(\prod_B C_i\right)c_n\in R(\s),$$ in which case we set
$$T_w=\left(\prod_{i\in B} T_i^V\right)T^V_\tau.$$ We then see that for $C_B,C_{B'}\in R(\s),$
we have
$$T_{C_B}T_{C_{B'}}=\prod_BT_i^V\prod_{B'}T_j^V=\prod_{B\wedge B'}T_i^V,$$
where $B\wedge B'$ is the symmetric difference.
Since $C_BC_{B'}=C_{B\wedge B'},$ we have the result in this case.  
A similar argument shows, in the case where $\G=D_n$ and $c_n\tau\simeq\tau,$ that  
$$T_{\bar C_B}T_{C_{B'}}=T_{\bar C_{B\wedge B'}}=T_{\bar C_BC_{B'}},$$
and
$$T_{\bar C_B}T_{\bar C_{B]}}=T_{C_{B\wedge B'}}=T_{\bar C_B\bar C_{B'}}.$$
Thus, we have the claim.\end{proof}

Since the cocycle $\eta:R(\s)\times R(\s)\rightarrow \C$ is determined by $T_{w_1w_2}=\eta(w_1,w_2)T_{w_1}T_{w_2}$ we have $\eta$ is a coboundary, and immediately get the following result.

\begin{corollary} For any Levi subgroup $\M$  of $\G_n$ and any $\s\in\mathcal E_2(M),$ we have $\mathcal C(\s)\simeq\C[R(\s)],$ so $i_{G,M}(\s)$ decomposes with multiplicity one.
\end{corollary}

Now, let $\A=\A_\theta$ be the split component of $\M,$ and let $\fa=\fa_\theta$ be its real Lie algebra.  If $\s\in\mathcal E_2(M),$ and $w\in R(\s),$ then we let $\fa_w=\{H\in\fa|w\cdot H=H\}.$  We let $Z$ be the split component of $G$ and $\fz$ be its real Lie algebra. Now, by Theorem 1.1 of \cite{herb},  we know $i_{G,M}(\s)$  has elliptic components if and only if there is a $w\in R(\s)$ with $\fa_w=\fz.$  Further, if
$\ds{\fa_{R(\s)}=\cap_{w\in R(\s)}}\fa_w,$ then each component of $i_{G,M}(\s)$ is irreducibly induced from an elliptic tempered representation if  there is some $w\in R(\s)$ so that $\fa_R=\fa_w.$

\begin{theorem}\label{elliptic odd} Let $\G=GSpin_{2n+1},$ and suppose $\M\simeq GL_{n_1}\times\cdots\times GL_{n_r}\times \G_m,$ and
$\s\in\mathcal E_2(M).$  Then $\Ind_P^G(\s)$ has elliptic constituents if and only if $R(\s)\simeq \Z_2^r.$  Any $\pi\in\mathcal E_t(G),$ is either elliptic, or there is a choice of $\M'$ and an irreducible elliptic tempered representation $\s$ of $M'$ with $\pi=\Ind_{P'}^G(\s).$
\end{theorem}

\begin{proof}We will use the explicit realization of $R(\s)$ we developed in Theorem ~\ref{R-gp odd}. Suppose $R\simeq\Z_2^d.$   
Let $\fa=\fa_M.$  Then we can identify $\fa$ with $\{(x_1,x_2,\dots,x_r,y)|x_i,y)\in \R \},$ and note, under this identificaton
$\fz=\{(y/2,\dots,y/2,y)|y\in \R\}.$
$\mathcal C$ acts on $\fa$ by
$$C_i(x_1,\dots x_{i-1},x_i,x_{i+1},\dots,x_r,y)=(x_1,\dots,x_{i-1},y-x_i,x_{i+1},\dots,x_r,y).$$
Thus, if $C=C_B,$ as above, then $\fa_C=\{(x_1,\dots,x_r,y)|x_i=y/2,\forall i\in B\}$
Without loss of generality, we may assume
$R(\s)=\langle C_r,C_{r-1},\dots, C_{r-d+1}\rangle.$  Let $w_0=C_{r-d+1}C_{r-d+2}\cdots C_r.$  Note, for each $w\in R(\s),$ we have $\fa_{w_0}\subset\fa_w,$ and thus $\fa_{R(\s)}=\fa_{w_0}.$  Now, $\fa_{w_0}=\fz$ if and only if $w_0=C_1C_2\cdots C_r,$ and thus, by \cite{arthur-elliptic,herb} $\Ind_P^G(\s)$ has elliptic constituents if and only if $R(\s)\simeq\Z_2^r.$  In this case, every component of $\Ind_P^G(\s)$ is elliptic.  The last statement of the claim follows from the fact $\fa_{R(\s)}=\fa_{w_0},$ and Lemma 1.3 of \cite{herb}.
\end{proof}

\begin{theorem}\label{elliptic even}
Let $\G=GSpin_{2n},$ and $\M\simeq GL_{n_1}\times\dots\times GL_{n_r}\times \G_m.$  Let $\s=\s_1\otimes\dots\otimes\s_r\otimes\tau\in\mathcal E_2(M).$
\begin{enumerate}[(i)]
\item Suppose $m=0$ or $c_n\tau\not\simeq\tau.$ We let $\O_1(\s),\O_2(\s),d_1,d_2,$ and $d$ be defined as in Theorem ~\ref{R-gp even}.  Then $\Ind_P^G(\s)$ has elliptic components if and only if  $d=r$ and $d_2$ is even, in which case every component is elliptic.  If $\pi\subset\Ind_P^G(\s)$ is not elliptic, then $\pi\simeq\Ind_{P'}^G(\s')$ for some elliptic representation $\s'$ of a Levi subgroup $M'$ of $G$ if and only if $d_2$ is even or $d_2=1.$
\item Suppose $m>0$ and $c_n\tau\simeq\tau.$  Let $R(\s)\simeq\Z_2^d.$  Then $\Ind_P^G(\s)$ has elliptic components if and only if  $d=r,$ in which case all components are elliptic.  Furthermore, for any $\pi\in\mathcal E_t(G)$ there is a Levi subgroup $\M'$ of $\G,$ and an irreducible elliptic tempered representation $\s'$ of $M'$ so that $\pi\simeq\Ind_{P'}^G(\s').$
\end{enumerate}
\end{theorem}

\begin{proof}  
\begin{enumerate}[(i)]
\item As in Theorem ~\ref{R-gp even} we assume $\O_1(\s)=\{r-d_1+1,r-d_1+2,\dots,r\},$ and $\O_2(\s)=\{r-d+1,r-d+2,\dots,r-d_1\}.$  Then
$$R(\s)=\left\langle C_iC_j|i,j\in\O_2(\s)\right\rangle\times\left\langle C_i|i\in\O_1(\s)\right\rangle.$$
We note $\fa=\fa_M$ can be identified with $\left\{\left(x_1,\dots,x_r,y\right) |x_i,y\in \R\right\},$ in such a way so
$C_i\cdot(x_1,\dots,x_{i-1},x_i,x_{i+1},\dots,x_r,y)=(x_1,\dots,x_{i-1},y-x_i,x_{i+1},\dots,x_r,y).$  So, if $d_2\neq1,$ we have
$$\fa_{R(\s)}=\left\{(x_1,\dots,x_r,y)|x_i=\frac{y}{2}\text{ for all } r-d+1\leq i\leq r\right\},$$
while if $d_2=1,$ then
$$\fa_{R(\s)}=\left\{(x_1,\dots,x_r,y)|x_i=\frac{y}{2}\text{ for all } r-d_1+1\leq i\leq r\right\}.$$
If $d_2$ is even then $w_0=C_{r-d+1}C_{r-d+2}\cdots C_r\in R(\s),$  and $\fa_{w_0}=\fa_{R(\s)}.$  If $d_2=1,$
then, $w_0=C_{r-d_1+1}C_{r-d+2}\cdots C_r\in R(\s),$ and again $\fa_{w_0}=\fa_{R(\s)}.$  
Thus, in either of these cases, we have each component must be irreducibly induced from an elliptic tempered representation of some Levi subgroup \cite{herb}.  On the other hand, if $d_2\geq 3$ and $d_2$ is odd, then, for any $w\in R(\s)$ we 
have $\fa_{R(\s)}\subsetneq\fa_w,$ so components of these induced representations are not irreducibly induced from elliptic representations.  Finally, since $\fz$ is identified with 
$\{(y/2,y/2,\dots,y/2,y)|y\in\R \},$ 
then we see $\Ind_P^G(\s)$ has elliptic components if and only if
$C_1C_2\dots C_r\in R(\s),$ which occurs if and only if $d=r$ and $d_2$ is even.

\item
Now suppose $m>0$ and $c_n\tau\simeq\tau.$   We let $\O(\s)$ be defined as in Theorem ~\ref{R-gp even}. We assume, without loss of generality,
$\O(\s)=\{r-d+1,\dots,r\}.$
Then
$$R(\s)=\langle C_i| r-d+1\leq i\leq r\text{ and } n_i \text{ is even}\rangle\times\langle C_ic_n| r-d+1\leq i\leq r\text{ and } n_i\text{ is odd}\rangle.$$
Let $d_2=\{i|r-d+1\leq i\leq r\text{ and }  n_i \text{ is odd}\},$  and
$$w_0=
\begin{cases} C_{r-d+1}C_{r-d+2}\cdots C_r&\text{ if } n_i\text{ is even};\\C_{r-d+1}C_{r-d+2}\cdots C_rc_n&\text{ if } n_i\text{ is odd}.
\end{cases}
$$
With the identification of $\fa=\fa_M$ with $\R^{r+1}$ as above, we have
$$\fa_{w_0}=\left\{\left(x_1\dots,x_r,y\right)|x_i=y/2\text{ for all } r-d+1\leq i\leq r\right\}.$$
Note, for any $w\in R(\s)$ we have $\fa_{w_0}\subset\fa_{w},$ so $\fa_{R(\s)}=\fa_{w_0}.$  Now, $\fa_{w_0}=\fz$ if and only if $d=r.$  Thus, the elliptic spectrum is as claimed, and the tempered spectrum is irreducibly induced from the elliptic spectra of the Levi subgroups.
\end{enumerate}
\end{proof}

 Now we assume $\G=\G_n=GSpin_{2n}$ or $GSpin_{2n+1}.$ Denote $R=R(\s),$ and let $\hat R$ be the set of  irreducible characters of $R.$ We let $\kappa\leftrightarrow\pi_\kappa$ be the correspondence between $\hat R$ and the (equivalence classes of)   irreducible components of $\Ind_P^G(\s)$ described by Keys \cite{keys-ens} (see also Arthur \cite{arthur-elliptic} and Herb \cite{herb}).  Suppose $\Ind_P^G(\s)$ has elliptic components, as described in Theorems ~\ref{elliptic odd}  and ~\ref{elliptic even}.  Then either $C_1C_2\dots C_r\in R$ or $C_1C_2\dots C_rc_n\in R.$ 
Let $$C_0=\begin{cases} C_1C_2\dots C_rc_n,&\text{ if } \G=GSpin_{2n}, d_2 \text{ is odd , and } c_n\tau\simeq\tau;\\ C_1C_2\dots C_n,&\text{ otherwise.}\end{cases}$$ 
For $\k\in\hat R$ we let $\e(\k)=\k(C_0).$
\begin{theorem}\label{elliptic character}
Suppose $\G=GSpin_{2n}$ or $GSpin_{2n+1}.$  Let $\M\simeq GL_{n_1}\times \dots\times GL_{n_r}\times \G_m$ be a Levi subgroup and suppose $\s=\s_1\otimes\dots\s_r\otimes\tau\in\mathcal E_2(M).$  Suppose $\Ind_P^G(\s)$ has elliptic components.   Let $\k\in\hat R.$  Then $\T^e_{\pi_\k}=\k(C_0)\T^e_{\pi_1}.$
\end{theorem}

\begin{proof}     First suppose $\G=GSpin_{2n+1},$ or $c_n\tau\simeq\tau.$ For $1\leq i\leq r,$ we let $\M_i$ be the Levi subgroup of $\G$ of the form $GL_{n_i}\times\G_{n-n_i}.$  Let $\N_i=\M_i\cap \N,$ and $\P_i=\M\N_i.$  We let $R_i=R_i(\s)$ be the $R$--group attached to $\Ind_{P_i}^{M_i}(\s).$  Since we are assuming $\Ind_P^G(\s)$ has elliptic components, we know $\Delta'=\emptyset.$  Thus, the compatibility condition in Section 2 of \cite{arthur-elliptic}  is satisfied (see also \cite{herb}).  Thus, we can identify  $R_i$ with the subgroup of $R$  generated by $\{C_j|1\leq j\leq r, j\neq i\},$ or 
$\{\bar C_j|1\leq j\leq r, j\neq i\},$  where $\bar C_i$ is define as in the proof of Theorem ~\ref{R-gp even}.  We now combine these situations by letting $R=\langle D_i|1\leq i\leq r\rangle,$ where $D_i=C_i$ or $\bar C_i,$ in the obvious way.
Let $\eta\leftrightarrow\rho_\eta$ be the correspondence between $\hat R_i$ and components of $\Ind_{P_i}^{M_i}(\s).$
If $\eta\in\hat R_i,$ we let $\hat R(\eta)=\{\k\in\hat R|\k|_{R_i}=\eta\}.$ Then $\hat R(\eta)=\{\eta^+,\eta^-\},$
where $\eta^\pm(D_j)=\eta(D_j),$ for $i\neq j,$ and $\eta^\pm(D_i)=\pm1.$  By Arthur \cite{arthur-elliptic} we have $\Ind_{M_iN_i'}^G(\rho_\eta)=\pi_{\eta^+}\oplus\pi_{\eta^-}.$  Moreover, since the character of this induced representation vanishes on $G_e,$ we have $\T_{\pi_{\eta^-}}^e=-\T^e_{\pi_{\eta^+}}.$ 

For $\k\in\hat R,$ we let $s(\k)=|\{1\leq i\leq r|\k(D_i)=-1\}|.$  Note, if $s(\k)=0,$ then $\k=1,$ and the claim is trivially true.  Suppose $s\geq 0$ and the claim holds for any $\k\in\hat R$ with $s(\k)=s.$  Suppose $s(\k)=s+1.$  Then we fix some $1\leq i\leq r$ for which $\k(D_i)=-1.$  Then consider $M_i$ and $R_i$ as above.  
Let $\eta=\k|_{R_i},$ and suppose  $\rho_\eta$  is the corresponding component of $\Ind_{P_i}^{M_i}(\s).$ Then $\k=\eta^-,$ so by our discussion above, we have $\T^e_{\pi_{\k}}=-\T^e_{\pi_{\eta^+}}.$  Moreover $s(\eta^+)=s,$ so, by our hypothesis, $\T^e(\pi_{\eta^+})=\eta^+(C_0)\T_1^e.$
Now, $\T_{\pi_\k}^e=-\T_{\pi_{\eta^+}}^e=-\eta^+(C_0)\T_1^e=\k(C_0)\T_e^1.$  So the claim holds for all $\k$ with $s(\k)=s+1,$ and by induction the claim holds for all $\k\in\hat R.$

Now consider the case where $\G=GSpin_{2n}$ and $c_n\tau\not\simeq\tau.$ The proof is essentially the same  as above, but we give some details for completeness.  Let $\O_1(\s),\,\O_2(\s) ,\, d_1,\,d_2,\, d$ be as in Theorem ~\ref{R-gp even}(i). If $d_2=0,$ then the proof is identical to the one above, so we assume $d_2>0$ is even.  From Theorem ~\ref{elliptic even}, we know $d=r.$  Then, we again see $\Delta'=\emptyset,$ so we easily apply the results of Arthur \cite{arthur-elliptic} and Herb \cite{herb}.  Without loss of generality, we assume $\O_1(\s)=\{1,\dots,d_1\},$ and $\O_2(\s)=\{d_1+1,\dots, r\}.$ Then, $R\simeq\Z_2^{r-1},$ with generators $D_1,\dots, D_{r-1},$ where $D_i=C_i,$ for $1\leq i\leq d_1,$ and $D_i=C_iC_r$ for $d_1+1\leq i\leq r-1.$
For each $1\leq i\leq r-1,$ we let $\M_i$ and $R_i$ be defined as in the previous cases.  We again let $\eta\leftrightarrow\rho_\eta$ be the correspondence between $\hat R_i$ and the components of $\Ind_{MN_i'}^{M_i}(\s).$  Then, we again have $\hat R(\eta)=\{\eta^+,\eta^-\},$ and so $\T_{\pi_{\eta^-}}^e=-\T_{\pi_{\eta^+}}^e.$  Let $\k\in\hat R$ and let $s(\k)=|\{D_i|\k(D_i)=-1\}|.$  Then $s(1)=0,$ so the claim holds for the case with $s(\k)=0.$  If we assume the result when $s(\k)=s,$ then the same argument as above  shows it holds when $s(\k)=s+1,$ and so the claim holds by induction.
\end{proof}
 
\

\section{Parameters and $R$-groups for $GSpin$ groups}

In this section we discuss the computation of Arthur's $R$-group associated to a parameter $\vp:W_G\rightarrow\,^LG,$ in the  case when $G=GSpin_m(F).$   
We begin with a lemma which applies to split reductive groups in general.

\begin{lemma}\label{maxcase1} Suppose $R_{\psi,\pi}\simeq R(\pi),$ whenever $\psi: W_F'\rightarrow\,^LL\hookrightarrow \,^LH,$ with $\bold L$ a maximal proper Levi subgroup of a quasi-split connected group $\bold H,$ and $\psi$ is an elliptic parameter for the $ L$-packet $\Pi_{\psi}(L),$ containing the square integrable representation $\pi.$  Let $\M$ be an arbitrary Levi subgroup of $\G,$ and $\vp:W_F'\rightarrow\,^LM$ an elliptic parameter for an $L$--packet $\Pi_{\vp}(M)$ containing a square integrable representation $\s.$ Then $R_{\vp,\s}\simeq R(\s).$
\end{lemma}

\begin{proof}

The proof of this relies on the following result.
\begin{lemma}\label{maxcase2}
Suppose $\M\subset \bold L$ are Levi subgroups of $\G.$  Suppose $\vp:W_F'\rightarrow\,^LM$ is a parameter.
Let $S_{\vp}=Z_{\hat G}(\vp)$ and $S_{L,\vp}=Z_{\hat L}(\vp).$  Then 
$S_{L,\vp}^\circ= S_{\vp}^\circ\cap\hat L.$
\end{lemma}
Since $S_\vp$ is reductive and $S_{L,\vp}$ is a reductive (Levi subgroup (e.g., by \cite{B-Zh} Lemma 2.1)
this is a standard result. 
\qed


Now we have $W(\G,\A_{\M})\simeq W(\hat G, A_{\hat M}),$ with the isomorphism given by $s_{\alpha}\mapsto s_{\ca}.$
We let $\M_{\alpha}$ be the Levi subgroup of $\G$ generated by $\M$ and $\alpha.$
Let $R_{\alpha}(\s)$ be the $R$-group attached to $i_{M_\alpha,M}(\s).$  
Considering $\vp:W_F'\rightarrow \,^LM\hookrightarrow\,^LM_{\alpha},$ we
 let $S_{\vp,\alpha}=Z_{\hat M_{\alpha}}(\vp)=S_\vp\cap \hat M_{\alpha}.$
By Lemma \ref{maxcase2}, $S_{\vp,\alpha}^\circ=S_{\vp}^\circ\cap\hat M_{\alpha}.$

We know from Lemma 2.2 of \cite{B-Zh} that $( A_{\hat M} \cap S_\vp)^\circ$ is a maximal torus 
of $ S_\vp^\circ$, so we may take $T_\vp = ( A_{\hat M} \cap S_\vp)^\circ$.
Then 
$^LM = Z_{^LG}(T_\vp)$ (\cite{B-Zh}, Lemma 2.1).
Since $T_\vp \subseteq \hat M \subseteq \hat M_{\alpha}$, it follows 
$T_\vp \subseteq S_{\vp,\alpha}$, so $T_\vp$ is a maximal torus in $S_{\vp,\alpha}^\circ$.

Let $W_{\vp,\alpha}=N_{S_{\vp,\alpha}}(T_{\vp})/Z_{S_{\vp,\alpha}}(T_{\vp}),$ and
$W_{\vp,\alpha}^\circ=N_{S_{\vp,\alpha}^\circ}(T_\vp)/Z_{S_{\vp,\alpha}^\circ}(T_\vp).$
Lemma 2.2 of \cite{B-Zh} tells us that $W_\vp$ (respectively, $W_{\vp,\alpha}$)
can be identified with the subgroup of $W(\hat G, A_{\hat M})$
(respectively, $W(\hat M_{\alpha}, A_{\hat M})$) consisting of the elements that
can be represented by elements of $S_\vp$ (respectively, $S_{\vp,\alpha}$).
Under these identifications, we have 
$
W_{\vp,\s}\cap \hat M_\alpha=W_{\vp,\alpha,\s}.$


Now let $R_{\vp,\alpha,\s}=W_{\vp,\alpha,\s}/W_{\vp,\alpha,\s}^\circ.$  The hypothesis implies $R_{\alpha}(\s)\simeq R_{\vp,\alpha,\s}.$
Let $\alpha\in\Delta'.$  Then $\mu_\alpha(\s)=0.$  Thus, $s_\alpha\in W(\s),$  and $R_{\alpha}(\s)=1.$
Note $s_{\alpha}\in W(\M_{\alpha},\A_{\M})\simeq W(\hat M_{\alpha}, A_{\hat M}), $ 
so $s_{\ca}\in W_{\vp,\s}\cap \hat M_\alpha=W_{\vp,\alpha,\s}.$  
Since $R_{\vp,\alpha,\s}\simeq R_{\alpha}(\s)=1,$ we have 
$s_{\ca}\in W_{\vp,\alpha,\s}^\circ,$  as claimed.
Conversely, assume $s_{\ca}\in W_{\vp,\s}^\circ.$  As $s_{\ca}\in W_{\vp,\s},$ we have $s_\alpha\in W(\s).$
Again, considering $\M_\alpha,$ we have $R_{\vp,\alpha,\s}=1,$ so $R_{\alpha}(\s)=1,$ which implies $s_\alpha\in W'.$ Therefore, $\alpha\in \Delta',$ as claimed.

\end{proof}

We now return to the setting where $\G=GSpin_m.$
Then  $\hat G=GSO_{2n}(\C)$, if $m=2n$, and
$\hat G=GSp_{2n}(\C)$, if $m=2n+1$.
Since $\G$ is split, we have $^LG=\hat G\times W_F.$  We consider a parameter 
$\vp:W_F\rightarrow\,^LG.$    Let us describe matrix realizations of $GSO_{2n}(\C)$ and
$GSp_{2n}(\C)$. Let
\[
  \mu = \begin{cases}
           1, & \text{if} \,  \hat G=GSO_{2n}(\C),\\
          -1, & \text{if} \,  \hat G=GSp_{2n}(\C),
        \end{cases}
\quad
 \hat w_n=  \small \begin{pmatrix} &&&&1\\&&&.\\&&.\\&.\\1\end{pmatrix}, \normalfont
\quad
  J_{2n} = \begin{pmatrix} 0&\hat w_n\\ \mu \hat w_n &0 \end{pmatrix},
\]
and 
$$
  \mathcal{G}= \{g\in GL_{2n}(\C)| \,^t g J_{2n}g=\l(g) J_{2n},\text{ for some }\l(g) \in\C^\times\}.
$$
If $\mu=-1$, then $\mathcal{G}$ is a connected algebraic group denoted by $GSp_{2n}(\C)$.
If $\mu=1$, then $\mathcal{G}=GO_{2n}(\C)$ has two connected components.
In this case, we can define the similitude norm
\[
    \nu : GO_{2n}(\C) \to \{ \pm 1 \}, \quad g \mapsto \lambda(g)^{-n} \det(g).
\]
The kernel of this map, denoted by $GSO_{2n}(\C)$, is the connected component of $GO_{2n}(\C)$.

We let $\hat M$ be the Siegel parabolic subgroup of $\hat G,$ so 
$\hat M\simeq GL_n(\C)\times GL_1(\C).$  More precisely, 
for $g\in GL_n(\C)$ we let $\hat\e(g)=\hat w_n\,^tg^{-1}\hat w_n^{-1}.$  Then
$$
   \hat M=\left\{\begin{pmatrix} g&0\\0&\l\hat\e(g)\end{pmatrix}
    \Bigg|g\in GL_n(\C),\l\in GL_1(\C)\right\}.
$$
Let $\hat A_{\hat M}$ be the split component of $\hat M,$ so $\hat M=\left\{\diag\left\{a I_n,\l a^{-1} I_n\right\}\right\}.$ 
If $\hat G=GSO_{2n},$ and  $n$ is odd, then $W_{\hat M}=\{1\}.$  Otherwise, 
$\hat W_{\hat M} =W(\hat G,\hat A_{\hat M})=\{1,\hat w\},$ where $\hat w:(g,\l)\mapsto (\l\hat\e(g),\l),$ and is represented by 
$
     \begin{pmatrix} 0& I_n\\ I_n&0\end{pmatrix}. 
$

Thus, we know $\M\simeq GL_n\times GL_1,\, \A_\M\simeq GL_1\times GL_1,$ and 
$$
W(\G,\A_{\M})=\begin{cases}
\{1\}&\text{if } G=GSpin_{2n}\text{~and } n\text{~is odd};\\
\{1,w\}&\text{otherwise,}
\end{cases}
$$
where $w:(g,\l)\mapsto (\l\e(g),\l),$ and $\e$ is the dual involution given by $\hat\e.$

Now let $\s$ be an irreducible unitary supercuspidal representation of $M,$ so 
$\s\simeq\s_0\otimes\psi,$ with $\s_0$ an irreducible unitary supercuspidal representation of 
$GL_n(F)$ and $\psi$ a unitary character of $F^\times.$  So, if $\vp:W_F\rightarrow\hat M$ is the corresponding 
Langlands parameter, then $\vp=\vp_0\times\hat\psi,$ where $\vp_0$ is the Langlands parameter of $\s_0$ and 
$\hat\psi$ is the character of $\C^\times$ associated to $\psi$ by local class field theory. 
Since $\s_0$ is irreducible and 
supercuspidal, we know $\vp_0$ is irreducible. We abuse notation to write
$$\vp(w)=\begin{pmatrix}\vp_0(w)&0\\0&\hat\psi(w)\hat\e(\vp_0(w))\end{pmatrix}.$$

\subsection{Reducibility and poles of $L$-functions}
Let $\hat{\mathfrak{n}}$ denote the Lie algebra of the unipotent radical of $\hat{M}$.
Let $\rho_n$ denote the standard representation of $GL_n(\mathbb{C})$.
The adjoint action $r$ of $\hat{M}$ on $\hat{\mathfrak{n}}$ is given as follows:
\[
  r = \begin{cases}
           \wedge^2 \rho_n \otimes \rho_1^{-1}, & \text{if} \, \hat G=GSO_{2n}(\C),\\
          \Sym^2 \rho_n \otimes \rho_1^{-1}, & \text{if} \, \hat G=GSp_{2n}(\C).
        \end{cases}
\]
More precisely, let $V = \{ X \in \mathfrak{gl}_{2n} (\mathbb{C}) \mid \, ^tX= - \mu X \}.$
Then $(g, \lambda) \in \hat M$ acts on $X \in V$ by 
$
   (g, \lambda) \cdot X = \lambda^{-1} gX \,^tg.
$

Suppose $L(s, \wedge^2 \vp_0 \otimes \hat\psi^{-1})$ has a pole at $s=0$. 
Then $ \wedge^2 \vp_0 \otimes \hat\psi^{-1}$ contains the trivial 
representation, so there exists a nonzero $X \in M_n(\C)$ such that $^tX = -X$ and
$ (\wedge^2 \vp_0 \otimes \hat\psi^{-1})(w) \cdot X = X$, for all $w \in W_F$. We have
\begin{equation}\label{formula}
   \hat\psi(w)^{-1} \vp_0(w) X \,^t\vp_0(w) = X, \quad \forall w \in W_F.
\end{equation}
It follows that $X$ is a nonzero intertwining operator between $^t\vp_0^{-1}$ and
$\hat\psi^{-1} \otimes \vp_0$. Since $\vp_0$ is irreducible, $X$ is invertible.
Observe that this can happen only if $n$ is even (every antisymmetric odd dimensional
matrix is singular).
In addition, it follows from (\ref{formula}) that $\vp_0$ factors through $GSp_n(\C)$.

Similarly, if we assume that $L(s, \Sym^2 \vp_0 \otimes \hat\psi^{-1})$ has a pole at $s=0$, we obtain that
$^t\vp_0^{-1} \simeq \hat\psi^{-1} \otimes \vp_0$ and $\vp_0$ factors through $GO_n(\C)$.

On the other hand, if $^t\vp_0^{-1} \simeq \hat\psi^{-1} \otimes \vp_0$, then 
(\ref{formula}) holds for some $X \in GL_n(\C)$. By standard arguments, $X$ is symmetric
or antisymmetric. It follows that one of the $L$-functions 
$L(s, \wedge^2 \vp_0 \otimes \hat\psi^{-1})$ or 
$L(s, \Sym^2 \vp_0 \otimes \hat\psi^{-1})$ has a pole at $s=0$.

We summarize the above considerations in the following lemma:

\begin{lemma}\label{factor}
Let $\vp_0 : W_F \to GL_n(\C)$ and $ \hat\psi : W_F \to GL_1(\C)$ be irreducible $L$-parameters.
If $\tilde{\vp}_0 \simeq \hat\psi^{-1} \otimes \vp_0$, then precisely one of the $L$-functions 
$L(s, \wedge^2 \vp_0 \otimes \hat\psi^{-1})$ or 
$L(s, \Sym^2 \vp_0 \otimes \hat\psi^{-1})$ has a pole at $s=0$.

\begin{enumerate}
  \item If $n$ is odd, then $L(s, \wedge^2 \vp_0 \otimes \hat\psi^{-1})$ is always holomorphic
   at $s=0$ and $\vp_0$ factors through $GO_n(\C)$.

\item If $n$ is even, then $L(s, \wedge^2 \vp_0 \otimes \hat\psi^{-1})$ has a pole
   at $s=0$ if and only if $\vp_0$ factors through $GSp_n(\C)$.
\end{enumerate}
\end{lemma}

\begin{proposition}\label{PropReduc}  

Let $\G=GSpin_{2n}$, $\G'=GSpin_{2n+1}$, and consider the Siegel Levi subgroup 
$\M\simeq GL_n\times GL_1.$  Let $\s\simeq\s_0\otimes\psi$ be an irreducible unitary supercuspidal representation of $M=\M(F)$ with corresponding Langlands parameter $\vp=\vp_0\times\hat\psi.$
Assume $\tilde{\vp}_0 \simeq \hat\psi^{-1} \otimes \vp_0$.
Let $\pi = \Ind_M^G(\s)$ and $\pi' = \Ind_M^{G'}(\s)$.

\begin{enumerate}
  \item If $n$ is odd, then $\pi$ and $\pi'$ are both irreducible and 
        $\vp_0$ factors through $GO_n(\C)$.
  
  \item If $n$ is even, then $\pi$  is irreducible if and only if $\pi'$ is reducible.
        Moreover, $\pi$  is irreducible if and only if $\vp_0$ factors through $GSp_n(\C)$.

\end{enumerate}
\end{proposition}

\begin{proof} (1) is clear. For (2), assume $n$ is even and consider $\G=GSpin_{2n}$ . Then 
$\pi = \Ind_M^G(\s)$ is irreducible if and only if 
$L(s, \s_0 \otimes\psi, \wedge^2 \rho_n \otimes \rho_1^{-1})$ has a pole at $s=0$
\cite{Sh2,silberger}. We know from \cite{Hen2}, Theorem 1.4 that
\[
   L(s, \s_0 \otimes\psi, \wedge^2 \rho_n \otimes \rho_1^{-1}) = 
     L(s, \wedge^2 \vp_0 \otimes \hat\psi^{-1} ).
\]
The statement follows from Lemma~\ref{factor}. Similar arguments work for $\G'=GSpin_{2n+1}$.
\end{proof}

\subsection{Centralizers for the Siegel Parabolic}

We wish to compute $S_{\vp} = Z_{\hat G}(\im\vp).$ 
First, we will compute $Z_{\mathcal G}(\im\vp)$, where ${\mathcal G} = GSp_{2n}(\C)$ or $GO_{2n}(\C)$.
 Suppose $X\in {\mathcal G}$ centralizes $\vp,$ and write $X=\begin{pmatrix} A&B\\C&D\end{pmatrix},$ with $A,B,C,D\in M_n(\C).$  Computing directly we have, for all $w\in W_F,$
$$
   \begin{pmatrix} A&B\\C&D\end{pmatrix}
      \begin{pmatrix} \vp_0(w)&0\\0&\hat\psi(w)\hat\e(\vp_0(w))\end{pmatrix} 
=
    \begin{pmatrix} \vp_0(w)&0\\0&\hat\psi(w)\hat\e(\vp_0(w))\end{pmatrix}
     \begin{pmatrix} A&B\\C&D\end{pmatrix},
$$ 
which gives
$$
A\vp_0(w)=\vp_0(w)A,\, D\hat\e(\vp_0(w))=\hat\e(\vp_0(w)) D,\, 
B\hat\psi(w)\hat\e(\vp_0(w))=\vp_0(w)B,\,
$$ 
and
$C\vp_0(w)=\hat\psi(w)\hat\e(\vp_0(w))C.$
The irreducibility of $\vp_0$ tells us $A$ and $D$ are scalars (denoted $a_{11} I_n$ and $a_{22} I_n,$ respectively) and also shows $C=B=0,$ unless $\vp_0\simeq(\hat\e\circ\vp_0)\otimes\hat\psi.$  
Thus, if $\s_0\not\simeq \tilde \s_0\otimes\psi\circ\det,$ 
then
$Z_{\mathcal G}(\vp)=\left\{\begin{pmatrix} a I_n\\&\l a^{-1}I_n\end{pmatrix}\right\}=\hat A_{\hat M}\simeq\C^\times\times\C^\times$, and
clearly, $Z_{\hat G}(\vp) = Z_{\mathcal G}(\vp)$.
So, suppose 
$\s_0\simeq \tilde\s_0\otimes\psi\circ\det.$  
Fix an equivalence, $B$ between these two representations, i.e.,
take $B$ so that $B^{-1}\vp_0(w)B=\hat\psi(w)\hat\e(\vp_0(w)).$  By Schur's Lemma, $B$ is unique up to scalar.  We note
$$(B\hat\e(B))^{-1}\vp_0(w)(B\hat\e(B))=\hat\e(B)^{-1}(\hat\psi(w)\hat\e(\vp_0)(w))\hat\e(B)=\hat\e(B^{-1}\vp_0(w)B)\hat\psi(w)=\vp_0(w),$$
and thus $B\hat\e(B)=cI_n,$ for some $c\in\C^\times.$  We write this as $B\hat w_n=c\hat w_n\,^tB.$  Note that if $J=B\hat w_n,$ then we have $^tJ=c^{-1}J,$ so $c=\pm 1,$ and $J$ is a symmetric or symplectic form fixed by $\vp_0$ up to the multiplier $\hat\psi.$

Now, we have
$X=\begin{pmatrix} a_{11}I_n&a_{12}B\\a_{21}B^{-1}&a_{22}I_n\end{pmatrix},$ and since 
$X\in \mathcal G,$ we have
\begin{gather*}
    ^t X\begin{pmatrix}&\hat w_n\\ \mu \hat w_n\end{pmatrix}X
            =\begin{pmatrix}&\l\hat w_n\\ \l \mu \hat w_n\end{pmatrix}\\
\intertext{or,}
\begin{pmatrix} a_{11}a_{21}(1+\mu c)\hat w_n B^{-1} & (a_{11}a_{22} + a_{21}a_{12}\mu c)\hat w_n
\\
             (a_{11}a_{22}+a_{12}a_{21} \mu c) \mu \hat w_n & a_{12}a_{22}(1+ \mu c)\,^tB \hat w_n
\end{pmatrix}
=
\begin{pmatrix}&\l\hat w_n\\ \l \mu\hat w_n\end{pmatrix}.
\end{gather*}
We see this is equivalent to the $2\times 2$ complex matrix $Y=\begin{pmatrix} a_{11}&a_{12}\\a_{21}&a_{22}\end{pmatrix}$ satisfying
$
^tY\begin{pmatrix} &1\\ \mu c \end{pmatrix} Y=\begin{pmatrix}&\l \\  \l \mu c \end{pmatrix}.
$  
Thus $X\mapsto Y$ is an isomorphism,
\begin{equation}\label{Zphi}
     Z_{\mathcal G}(\vp)\simeq  
   \begin{cases} GSp_{2}(\C)\simeq GL_2(\C)&\text{ if } \mu c=-1;\\
                      GO_{1,1}(\C)&\text{ if } \mu c=1.
   \end{cases}
\end{equation}
This is equal to $S_{\vp}$ if $\hat G = GSp_{2n}(\C)$.  

Now, let  $\hat G = GSO_{2n}(\C)$, so $\mu = 1$.
Let
$X=\begin{pmatrix} a_{11}I_n&a_{12}B\\a_{21}B^{-1}&a_{22}I_n\end{pmatrix} \in Z_{\mathcal{G}}(\vp)$.
We have to determine whether $X \in \hat G$.
Assume first $c= -1$. Then 
\[
\begin{pmatrix}  &(a_{11}a_{22}-a_{12}a_{21})\hat w_n
\\
  (a_{11}a_{22}-a_{12}a_{21})\hat w_n&  \end{pmatrix}
=
    \begin{pmatrix}&\l\hat w_n\\\l\hat w_n\end{pmatrix},
\]
so $\lambda = a_{11}a_{22}-a_{12}a_{21}$.
We use the formula
$
  \det \begin{pmatrix} A & B\\C & D \end{pmatrix}
=
  \det A \det(D- CA^{-1}B),
$
if $A$ is invertible. Therefore, if $a_{11} \ne 0$, we have
\[
   \det X = a_{11}^n \det( a_{22}I_n - a_{21}a_{11}^{-1}a_{12}B^{-1}B)
      = \det(a_{11}a_{22}-a_{12}a_{21})I_n
      = \lambda^n.
\]
The similitude norm $\nu(X) = \lambda^{-n} \det X = 1$, so $X \in GSO_{2n}(\C)$.
If $a_{11} = 0$, then
\[
   \det X = \det \begin{pmatrix} 0 &a_{12}B\\a_{21}B^{-1}&a_{22}I_n\end{pmatrix}
   = (-1)^n \det \begin{pmatrix} a_{21}B^{-1}&a_{22}I_n \\ 0 &a_{12}B \end{pmatrix}
      = \lambda^n,
\]
and again $X \in GSO_{2n}(\C)$.

Assume $c=1$. Then we have
\[
\begin{pmatrix} 2a_{11}a_{21} \hat w_n B^{-1} &(a_{21}a_{12}+a_{11}a_{22})\hat w_n
\\
  (a_{11}a_{22}+a_{12}a_{21})\hat w_n&  2a_{12}a_{22}\,^tB \hat w_n\end{pmatrix}
=
    \begin{pmatrix}&\l\hat w_n\\\l\hat w_n\end{pmatrix}.
\]
It follows $a_{12} = a_{21} =0$ or $a_{11} = a_{22} =0$. If $a_{12} = a_{21} =0$,
then $a_{22} = \lambda a_{11}^{-1}$ and 
$
X=\begin{pmatrix} a_{11}I_n & \\ & \lambda a_{11}^{-1}I_n \end{pmatrix}.
$
The similitude norm
$
   \nu(X) = \lambda^{-n} \det(X) = 1,
$
so $X \in GSO_{2n}(\C)$. If $a_{11} = a_{22} =0$, then 
$
 X=\begin{pmatrix}  & a_{12}B \\ \lambda a_{12}^{-1}B^{-1}& \end{pmatrix}
$
and
\[
   \nu(X) = \lambda^{-n} \det(X) = (-1)^n \lambda^{-n} \lambda^{n} = (-1)^n.
\]
It follows that $X \in GSO_{2n}(\C)$ if $n$ is even and $X \notin GSO_{2n}(\C)$ if $n$ is odd.
Therefore, 
 $$
S_\vp=Z_{\hat G}(\vp)\simeq
\begin{cases} GSp_{2}(\C)\simeq GL_2(\C)&\text{ if }c=-1;
\\
GO_{1,1}(\C)&\text{ if }c=1, \, $n$ \text{ even}; 
\\
\C^\times &\text{ if }c=1, \, $n$ \text{ odd}.
\end{cases}
$$

\subsection{The Arthur $R$-group}

Now we can compute $R_{\vp}$, the Arthur $R$-group of ${\vp}$. We summarize the above 
computation as follows.

\begin{theorem} Let $\G=GSpin_{2n+1}$ or $GSpin_{2n}$  and consider the Siegel Levi subgroup 
$\M\simeq GL_n\times GL_1.$  Let $\s\simeq\s_0\otimes\psi$ be an irreducible unitary supercuspidal representation of $M=\M(F)$ with corresponding Langlands parameter $\vp=\vp_0\otimes\hat\psi.$ 
\begin{enumerate}
  \item If $\vp_0\not\simeq \tilde \vp_0\otimes\hat\psi$, then $R_{\vp,\s}=R_{\vp}= 1$.
  \item Assume $\vp_0 \simeq \tilde \vp_0\otimes\hat\psi$. If $\G=GSpin_{2n+1}$, then
\[
   R_{\vp,\s} = R_{\vp}=\begin{cases}
          1, &  \text{ if }\vp_0\text{ factors through } GO_n(\C);\\
          \Z_2, & \text{ if } \vp_0\text{ factors through } GSp_n(\C).
    \end{cases}
\]
If $\G=GSpin_{2n}$, then
\[
   R_{\vp,\s} = R_{\vp}=\begin{cases}
          1, &  \text{ if }\vp_0\text{ factors through } GSp_n(\C);\\
          \Z_2, & \text{ if } \vp_0\text{ factors through } GO_n(\C) \text{ and } n \text{ is even,}\\
          1, & \text{ if } \vp_0\text{ factors through } GO_n(\C) \text{ and } n \text{ is odd.}    
\end{cases}
\]

\end{enumerate}

\end{theorem}

\begin{corollary}\label{corollary siegel}
For $\G=GSpin_{2n+1}$ or $GSpin_{2n}$, and $\M\simeq GL_n\times GL_1,$ we have
$R(\s)\simeq R_{\varphi,\s},$ as conjectured by Arthur.
\end{corollary}

\begin{proof}
This follows from the theorem and Proposition~\ref{PropReduc}.
\end{proof}

\subsection{Centralizers (The General Case)}
Let $V$ be a finite dimensional complex vector space.
Let $B$ be a non-degenerate bilinear form on $V$ and  
$$
  \mathcal{G}_B= \{g\in GL_{n}(V) \mid B(gu, gv) = \lambda(g)B(u,v) ,\text{ for some }\l(g) \in\C^\times, \, \forall u,v \in V \}.
$$

\begin{lemma}\label{bilin}

Let $\vp : W_F' \to GL_{n}(V) $ be an irreducible parameter and let $B$ be a non-degenerate bilinear form on $V$.
Then $\vp$ factors through $\mathcal{G}_B$ if and only if
 $\vp \simeq \chi \otimes \,^t\vp ^{-1}$, where $\chi = \lambda \circ \vp$.
If $\vp$ factors through $\mathcal{G}_B$, then $B$ is unique up to a scalar multiple.

\end{lemma}

\begin{proof} Suppose that $\vp$ factors through $\mathcal{G}_B$. Let $A$ be the matrix corresponding to $B$, $B(u,v) = \,^t u A v$. Then for all $w \in W_F'$,
$
     ^t \vp(w)A \vp(w) = \l(\vp(w))A.
$
It follows 
\begin{equation}\label{intertw}
     \vp(w) = \chi(w) A \,  ^t \vp(w)^{-1}A^{-1}, \quad \forall w \in W_F',
\end{equation}
 where $\chi = \lambda \circ \vp$. Hence, $\vp \simeq \chi \otimes \,^t\vp^{-1}$.
If $B'$ is another non-degenerate bilinear form on $V$
such that $\vp$ factors through $\mathcal{G}_{B'}$, and if $A'$ is the corresponding matrix,
we have 
\begin{equation}\label{intertw2}
     \vp(w) = \chi(w) A' \,  ^t \vp(w)^{-1}(A')^{-1}, \quad \forall w \in W_F'.
\end{equation}
By transposing and taking inverses, equation (\ref{intertw}) gives us
$
   ^t \vp(w)^{-1} = \chi(w)^{-1}A^{-1} \vp(w) A.
$
We substitute this in equation (\ref{intertw2}) and we obtain
\[
     \vp(w) = A' A^{-1} \vp(w) A (A')^{-1}, \quad \forall w \in W_F'.
\]
Since $\vp$ is irreducible, it follows $A' A^{-1}=cI$ and $A'=cA$.

Next, suppose $\vp \simeq \chi \otimes \,^t\vp ^{-1}$ for a character $\chi$.
Let $A$ be a matrix such that 
\[
     \vp(w) = \chi(w) A \,  ^t \vp(w)^{-1}A^{-1},
\]
for all $w \in W_F'$. Standard arguments show that $A \, ^tA^{-1} = cI$ and $c=\pm 1$.
It follows that $B(u,v) = \,^t u A v$ is a  non-degenerate bilinear form such that $\vp$ factors
through $\mathcal{G}_B$.
\end{proof}

\begin{lemma}\label{LemmaZ} 
Let $\vp : W_F' \to \mathcal{G}_B$ be a parameter. Suppose
$\vp = \underbrace{\vp_0 \oplus \cdots \oplus \vp_0}_{m-\text{summands}}$,
where $\vp_0$ is an irreducible parameter such that $\vp_0$ factors through $\mathcal{G}_{B_0}$
for some non-degenerate bilinear form $B_0$. Then
\[
     Z_{\mathcal{G}_B}(\im \vp) \simeq \begin{cases}
                GO(m, \C), & if \, B \, and \, B_0 \, \text{are both symmetric or both symplectic,}\\
                GSp(m,\C), & \text{otherwise.}
    \end{cases}
\]

\end{lemma}

\begin{proof} Let $V_0$ denote the space of the representation $\vp_0$.
Then $V \simeq W \otimes V_0$, where
$
 W = \operatorname{Hom}_{W_F'}(V_0,V)
$
with trivial $W_F'$-action. The map $W \otimes V_0 \to V$ is given by
\begin{equation}\label{ident}
  f \otimes v \mapsto f(v), \quad f \in W, v \in V_0.
\end{equation}
For $f,g \in W$, we define a bilinear form $B_{f,g}$ on $V_0$ by
$B_{f,g}(u,v) = B(f(u), g(v))$. Then
\[
\begin{aligned}
    B_{f,g}(\vp_0(w)u,\vp_0(w)v) &= B(f(\vp_0(w)u),g(\vp_0(w)v)) \\
&=
     B(\vp(w)f(u),\vp(w)g(v))\\
&= \lambda \circ \vp (w) B_{f,g}(u,v).
\end{aligned}
\]
It follows from Lemma~\ref{bilin} that $B_{f,g}$ is a scalar multiple of $B_0$;
denote that scalar by $\langle f, g \rangle$. The map $(f,g) \mapsto \langle f, g \rangle$
defines a bilinear form $\langle\, , \rangle$ on $W$. The form $\langle \, , \rangle$ is symmetric if 
$B$ and $B_0$ are both symmetric or both symplectic, and symplectic otherwise.
Moreover, if we identify $W \otimes V_0$ and $V$ using equation (\ref{ident}), we have
\[
   B(f \otimes u, g \otimes v) = B(f(u), g(v)) = B_{f,g}(u,v) = \langle f, g \rangle B_0(u,v),
\]
for all $f,g \in W$, $u,v \in V_0$.

Now, $\im \vp = \{ I_W \otimes g \mid g \in \im \vp_0 \}$ and
\[
   Z_{GL(V)}(\im \vp) = \{ g \otimes z \mid g \in GL(W), z = cI_{V_0}, c \in \C^\times \}
        = \{ g \otimes I_{V_0} \mid g \in GL(W) \}.
\]
The element $g \otimes I_{V_0}$ belongs to $\mathcal{G}_B$ if for some $\lambda \in \C^\times$,
\[
  B((g \otimes I_{V_0}) (f \otimes u), (g \otimes I_{V_0}) (h \otimes v))
     = \lambda B(f \otimes u, h \otimes v), \quad \forall f, h \in W, \forall u,v \in V_0,
\]
that is,
\[
    \langle g f, g h \rangle = \lambda \langle f, h \rangle, \quad \forall f, h \in W.
\]
It follows $Z_{\mathcal{G}_B}(\im \vp) \simeq \mathcal{G}_{\langle \,, \rangle}$.
\end{proof}

\subsection{Reducibility for generic representations}

Let $G=GSpin_m(F)$ and let $P=MN$ be a maximal Levi subgroup. Then 
$M \simeq GL_k(F) \times GSpin_\ell(F)$, where $2k + \ell = m$.
In the case $\ell = 0$ or 1, $P$ is the Siegel parabolic subgroup and that case was considered
earlier. We assume $\ell > 2$.
Let $\pi = \sigma \otimes \tau$ be an irreducible unitary generic supercuspidal
representation of $M$. Let $\alpha \in \Delta$ be the unique reduced root of $\T$ in $\N$
and set $\tilde{\alpha} = \langle \rho, \alpha \rangle^{-1} \alpha$, where $\rho$ is half
the sum of positive roots in $\N$. We have 
$ 
 \tilde{\alpha}/i \otimes \pi = \nu^{1/i} \sigma \otimes \tau.
$
Assume $\sigma \simeq \tilde{\sigma} \otimes \omega_\tau$.
According to \cite{Sh2}, exactly one of the following representations
is reducible: $\Ind_P^G(\sigma \otimes \tau)$, $\Ind_P^G(\nu^{1/2}\sigma \otimes \tau)$, 
or $\Ind_P^G(\nu\sigma \otimes \tau)$.

\begin{lemma}\label{factor2}

Let $G=GSpin_m(F)$ and $M \simeq GL_k(F) \times GSpin_\ell(F)$, where $2k + \ell = m$, $\ell > 2$.
Let $\pi = \sigma \otimes \tau$ be an irreducible unitary generic supercuspidal
representation of $M.$ Assume $\sigma \simeq \tilde{\sigma} \otimes \omega_\tau$.
Let $\vp_0$ denote the Langlands parameter of $\sigma$.

\begin{enumerate}
  \item Suppose $G=GSpin_{2n+1}(F)$. If $\Ind_P^G(\nu^{1/2}\sigma \otimes \tau)$ reduces, 
        then $\vp_0$ factors through $GO_k(\C)$. Otherwise,  $\vp_0$ factors through $GSp_k(\C)$.
  \item Suppose $G=GSpin_{2n}(F)$. If $\Ind_P^G(\nu^{1/2}\sigma \otimes \tau)$ reduces, 
        then $\vp_0$ factors through $GSp_k(\C)$. Otherwise,  $\vp_0$ factors through $GO_k(\C)$.
 
\end{enumerate} 

\end{lemma}

\begin{proof}

Let $\hat{\mathfrak{n}}$ denote the Lie algebra of the unipotent radical of $\hat{M}$.
Denote the standard representations of the groups $GL_k(\C)$, $GSp_{2\ell}(\C)$
and $GSO_{2\ell}(\C)$ by $\rho_k$, $R^1_{2\ell}$ and $R^2_{2\ell}$, respectively.
Let $\mu$ be the similitude character of $GSp_{2\ell}(\C)$ or $GSO_{2\ell}(\C)$.
The adjoint action $r$ of $\hat{M}$ on $\hat{\mathfrak{n}}$ is described
in Proposition 5.6 of \cite{asgari}.  In particular, we have

\begin{enumerate}
\item[{(a)}] If $G=GSpin_{2n+1}(F)$, then $r= r_1 \oplus r_2$, where
$$
    r_1= \rho_k \otimes \widetilde{R^1_{\ell-1}}, \quad 
    r_2 = \Sym^2 \rho_k \otimes \mu^{-1}.
$$
\item[{(b)}]  If $G=GSpin_{2n}(F)$, then $r= r_1 \oplus r_2$, where
$$
    r_1= \rho_k \otimes \widetilde{R^2_{\ell}}, \quad 
    r_2 = \wedge^2 \rho_k \otimes \mu^{-1}.
$$

\end{enumerate}

Let $P_{\pi,1}$ and $P_{\pi,2}$ be the polynomials defined in \cite{Sh2}. 
The Langlands-Shahidi $L$-function attached to $\pi$ and $r_i$ is defined as
\[
    L(s, \pi, r_i) = P_{\pi,i}(q^{-s})^{-1}.
\]
Assume $G=GSpin_{2n}(F)$.
Theorem 8.1 of \cite{Sh2} tells us that  $\Ind_P^G(\nu^{1/2}\sigma \otimes \tau)$ is reducible
if and only if $P_{\pi,2}(1)=0$. Equivalently, $ L(s, \pi, r_2) $ has a pole at $s=0$.
In order to complete the proof, we need the following result.
\begin{lemma}\label{2nd lf}

Let $\G=GSpin_{m}$ and $\M\simeq GL_k\times GSpin_{\ell}.$   Let $\pi = \sigma \otimes \tau$ be an irreducible admissible generic 
representation of $M.$  Let $\vp=(\vp_0, \vp_{\tau})$ be the Langlands parameter attached to $\pi.$

\begin{enumerate}[a)] 
\item If $m=2n$ is even, then
\begin{equation}\label{claim-even}
 \begin{aligned}
    L(s, \pi, r_2) &= L(s, \sigma \otimes \tau, \wedge^2 \rho_k \otimes \mu^{-1})\\
   &=L(s, \sigma \otimes \psi, \wedge^2 \rho_k \otimes \rho_1^{-1}) 
   = L(s, \wedge^2 \vp_0 \otimes \hat{\psi}^{-1}).
\end{aligned}
\end{equation}

\item If $m=2n+1$ is odd, then
\begin{equation}
 \begin{aligned}
    L(s, \pi, r_2) &= L(s, \sigma \otimes \tau, \Sym^2 \rho_k \otimes \mu^{-1})\\
   &=L(s, \sigma \otimes \psi, \Sym^2 \rho_k \otimes \rho_1^{-1}) 
   = L(s, \Sym^2 \vp_0 \otimes \hat{\psi}^{-1}).
\end{aligned}
\end{equation}

\end{enumerate}
\end{lemma}
\begin{proof}
 We continue with the notation of the proof of Lemma \ref{factor2}.
 Suppose $m=2n.$ Then $\vp_\tau:W_F\rightarrow GSO_{2\ell}(\C).$   
First, we prove \eqref{claim-even} holds for any unramified generic $\pi$.
By Prop. 2.3(a) of \cite{asgari-shahidi-djm} we know $Z(GSp_{2\ell}(F))^\circ=\{e_0^*(\l)|\l\in F^\times\}.$  So, the central character of $\tau$ is given by $$\o_\tau(\l)\operatorname{Id}_{V_\tau}=\tau(e_0^*(\l)).$$
Let $\hat\psi:W_F\rightarrow \C^\times$ be the character attached to $\o_\tau$ by Class Field Theory.  In particular, $\o_\tau(\varpi_F)=\hat\psi(\Fr_F),$
where $\Fr_F$ is the Frobenius class of $F.$
Let $\hat T$ be the maximal torus of $GSO_{2m}(\C).$  Then $\mu(t)=e_0^*(t),$ (by \cite{asgari-shahidi-djm}  pg. 149).    Now, we have
$$
L(s,\pi,r_2)=L(s,\s\otimes\tau,\wedge^2\rho_k\otimes\mu^{-1})
    =L(s,\wedge^2\rho_k\otimes\mu^{-1}(\vp_0,\vp_\tau)).
$$
Note, for  $w\in W_F,$ we have
 $\wedge^2\rho_k\otimes\mu^{-1}(\vp_0,\vp_\tau)(w)=\wedge^2(\vp_0(w))\mu^{-1}(\vp_\tau(w)).$
Now 
$$
   \mu^{-1}(\vp_\tau(\Fr_F))=(e_0^*(\vp_\tau(\Fr_F)))^{-1}
     =\tau(e_0^*(\varpi_F))^{-1}=\o_\tau(\varpi_F)^{-1}.
$$
So $$L(s,\pi,r_2)=L(s,\wedge^2\rho_k\vp_0\otimes\hat\psi^{-1})=L(s,\s\otimes\o_\tau,\wedge^2\rho_k\otimes\rho_1^{-1}).$$

If $S_n$ denotes the $n$-dimensional complex representation of $SL(2,\C)$,
then $\im(S_n)$ is orthogonal or symplectic. Therefore, 
$
   \mu ( \vp \otimes S_n) = \mu (\vp).
$
We conclude that equation \eqref{claim-even} holds if $\pi$ has an Iwahori fixed vector. In addition, for $\vp$ unramified,
the Artin $\varepsilon$--factor associated to $ \mu ( \vp \otimes S_n)$ is equal to 1.

Now, we apply Theorem 3.5 of \cite{Sh2} to $\pi=\sigma \otimes \tau$ and 
independently we apply the same theorem to $\sigma \otimes \o_\tau$.
The theorem guarantees existence of the $\gamma$--factors 
$\gamma_2(s, \sigma \otimes \tau, \psi_F, \tilde{w})$ and 
$\gamma_1(s, \sigma \otimes \o_\tau, \psi_F, \tilde{w}),$ with the subscripts determined by the ordering of the components of the adjoint representations of the $L$--groups of the Levi subgroups in two distinct situations.
Moreover, conditions 1, 3, and 4 from the theorem determine these $\gamma$--factors uniquely.
These conditions are satisfied  by $\gamma_2(s, \sigma \otimes \tau, \psi_F, \tilde{w})$ and 
independently by
$\gamma_1(s, \sigma \otimes \o_\tau, \psi_F, \tilde{w})$.
In the inductive property 3 for $\gamma_1(s, \sigma \otimes \o_\tau, \psi_F, \tilde{w})$,
only $\sigma$ can be induced from a smaller parabolic subgroup, not $\o_\tau$.
Therefore, if we look at the inductive property for 
$\gamma_1(s, \sigma \otimes \o_\tau, \psi_F, \tilde{w})$, the same conditions are 
satisfied for $\gamma_2(s, \sigma \otimes \tau, \psi_F, \tilde{w})$.
Even though we can have additional conditions for 
$\gamma_2(s, \sigma \otimes \tau, \psi_F, \tilde{w})$, the conditions for 
$\gamma_1(s, \sigma \otimes \o_\tau, \psi_F, \tilde{w})$ are enough to guarantee uniqueness.
Since we have equality of  $\gamma$--factors for representations with Iwahori fixed vectors,
we conclude that 
$\gamma_2(s, \sigma \otimes \tau, \psi_F, \tilde{w}) = 
\gamma_1(s, \sigma \otimes \o_\tau, \psi_F, \tilde{w})$.
The definition of $L$-functions from \cite{Sh2} then implies \eqref{claim-even}.

The proof of the case $\G=GSpin_{2n+1}$ is similar.
\end{proof}

We return to the proof of Lemma \ref{factor2}. It follows from Lemma~\ref{factor} and Lemma \ref{2nd lf}  that $\Ind_P^G(\nu^{1/2}\sigma \otimes \tau)$ is reducible
if and only if $\vp_0$ factors through $GSp_k(\C)$.
Finally, we remark the claim will follow in general from the generic $L$--packet conjecture of Shahidi \cite{Sh2}.

The proof for $G=GSpin_{2n+1}(F)$ is similar. 
\end{proof}

Let $G=GSpin_{2\ell+1}(F)$ and let $\tau$ be a generic discrete series
representation of $G$. As in \cite{Moe}, let $\Jord(\tau)$ denote the set of pairs $(\rho, a)$,
where $\rho \in \,^0\mathcal{E}(GL(d_\rho,F))$ and $a \in \Z^+$ such that 
$\delta(\rho,a) \rtimes \tau$ is irreducible and there exists an integer $a'$ of the same
parity as $a$ such that $\delta(\rho,a') \rtimes \tau$ is reducible.
Here $\delta\rtimes\tau=\Ind^{G(\ell+d)}_{GL_d(F)\times G(\ell)}(\delta\otimes\tau).$
The $L$-parameter of $\tau$ is given by
\[
    \vp_\tau = \bigoplus_{(\rho,a) \in \Jord(\tau)} \vp_\rho \otimes S_a,
\]
where $\vp_\rho$ is the $L$-parameter of $\rho$.

\begin{theorem}\label{isomorphism max-non-siegel}

Let $\G=GSpin_{2n+1}$  and consider the Levi subgroup 
$M \simeq GL_k(F) \times GSpin_{2\ell+1}(F)$.  
Let $\pi = \s \otimes \tau$ be a generic discrete series representation of $M$.
 Let $\vp $ be the $L$-parameter of $\pi$. Then $R_{\vp, \pi} \simeq R(\pi)$.

\end{theorem}

\begin{proof}

The parameter $\vp$ can be written as 
$\vp \simeq \vp_\s \oplus \vp_\tau \oplus (\hat\e({\vp}_\s) \otimes \hat{\psi})$,
where $\hat{\psi}$ is  the character corresponding to the central character of $\tau,$ (restricted to the connected component of the center)  by Class Field Theory.
The representation $\s$ is of the form $\s \simeq \delta(\rho, a)$, where
$\rho \in \,^0\mathcal{E}(GL(d,F))$ and $a \in \Z^+$, $da=k$. Then 
$\vp_\sigma = \vp_\rho \otimes S_a$.

If $\sigma \not \simeq \tilde{\sigma} \otimes \omega_\tau$, it is easy to show 
$R_{\vp, \pi}=1$ and $R(\pi)=1$.
Assume $\sigma \simeq \tilde{\sigma} \otimes \omega_\tau$. Then 
$\hat\e({\vp}_\s) \otimes \hat{\psi} \simeq \vp_\sigma$, so 
$\vp \simeq \vp_\s \oplus \vp_\tau \oplus \vp_\s$.

If $(\rho,a) \in \Jord(\tau)$, then the multiplicity of $\vp_\s$
in $\vp \simeq \vp_\s \oplus \vp_\tau \oplus \vp_\s$ is three.
Lemma~\ref{LemmaZ} implies $R_\vp =1$. On the other hand, since $(\rho,a) \in \Jord(\tau)$,
we have $\s \rtimes \tau$ is irreducible, so $R(\pi)=1$.

Now, consider the case $\sigma \simeq \tilde{\sigma} \otimes \omega_\tau$ and
$(\rho,a) \notin \Jord(\tau)$. There exist a supercuspidal generic representation $\tau_{cusp}$
of $GSpin_{2m+1}(F)$ and an irreducible generic representation $\theta$ of $GL_r(F)$ such that
$\tau$ is a subrepresentation of
$$
     \theta \rtimes \tau_{cusp} =i_{GSpin_{2\ell+1}(F), GL_r(F) \times GSpin_{2m+1}(F)}
     (\theta \otimes \tau_{cusp}).
$$
We apply the Langlands classification for $GL_r(F)$ in the subrepresentation setting. 
It follows that there exist 
$\delta(\rho_1, a_1), \delta(\rho_2, a_2), \dots, \delta(\rho_s, a_s)$
and real numbers $b_1 < b_2 < \cdots < b_s$ such that $\theta$ is the unique subrepresentation
of  the induced representation
$$
    \nu^{b_1} \delta(\rho_1, a_1)\times  \nu^{b_2} \delta(\rho_2, a_2) \times
     \cdots \times \nu^{b_s} \delta(\rho_s, a_s).
$$
For $i \in \{1, \dots ,s\}$, define $[i]=\{ j  \in \{1, \dots ,s\} \mid \rho_i \simeq \rho_j \}$.
The Casselman square integrability criterion for $\tau$ implies that for $i=1, \dots ,s$,
there exists $j \in [i]$ such that the representation
$$
        \nu^{b_j} \delta(\rho_j, a_j) \rtimes \tau_{cusp}
$$
is reducible, and $b_i - b_j \in \Z$. This implies $b_j \in \frac{1}{2}\Z$ and therefore
 $b_i \in \frac{1}{2}\Z$.

Assume first $\s \rtimes \tau$ is reducible. Then $R(\pi)\simeq \Z_2$.
It can be shown, taking into account the structure of $\theta$, 
that reducibility of $\s \rtimes \tau$ implies reducibility of $\s \rtimes \tau_{cusp}$.
Then there exists $b \geq 0$, 
$b \in \{-\frac{(a-1)}{2}, -\frac{(a-1)}{2} + 1, \dots , \frac{(a-1)}{2}\}$
such that $\nu^b \rho \rtimes \tau_{cusp}$ is reducible.
Since $\tau_{cusp}$ is supercuspidal 
and generic, we have $b=0, 1/2$ or 1. If $b=1/2$, then $a$ is even. In addition, 
Lemma~\ref{factor2} implies that $\vp_\rho$ factors through $GO_d(\C)$. Then 
$\vp_\sigma = \vp_\rho \otimes S_a$ factors through $GSp_k(\C)$. Now Lemma~\ref{LemmaZ} 
tells us that $S_\vp \simeq GO(2,\C)$. It follows $R_\vp = R_{\vp, \pi}\simeq \Z_2$.
If $b=0$ or 1, then $a$ is odd. In addition, 
Lemma~\ref{factor2} implies that $\vp_\rho$ factors through $GSp_d(\C)$. Then 
$\vp_\sigma = \vp_\rho \otimes S_a$ factors through $GSp_k(\C)$.
As before, we obtain $R_{\vp, \pi}\simeq \Z_2$.

It remains to consider the case when $\s \rtimes \tau$ is irreducible, 
$\sigma \simeq \tilde{\sigma} \otimes \omega_\tau$ and $(\rho,a) \notin \Jord(\tau)$.
Irreducibility of $\s \rtimes \tau$ implies $R(\pi)=1$.
Let $b \in \{0, 1/2, 1 \}$ such that $\nu^b \rho \rtimes \tau_{cusp}$ is reducible.
Since $(\rho,a) \notin \Jord(\tau)$, $a$ and $2b+1$ are not of the same parity.
Therefore, if $b=1/2$, then $a$ is odd. Then $\vp_\rho$ factors through $GO_d(\C)$ and 
$\vp_\sigma = \vp_\rho \otimes S_a$ factors through $GO_k(\C)$. It follows
$R_\vp = R_{\vp, \pi}=1$. Similarly, if $b=0$ or 1, then $a$ is even, 
$\vp_\rho$ factors through $GSp_d(\C)$ and 
$\vp_\sigma = \vp_\rho \otimes S_a$ factors through $GO_k(\C)$, implying 
$R_\vp = R_{\vp, \pi}=1$.
\end{proof}

\begin{theorem}\label{isomorphism general case}
Let $\G=GSpin_{2n+1}$ and $\P=\M\N$ be an arbitrary parabolic subgroup of $\G.$  Suppose $\pi$ is a discrete series representation of $M$ and $\vp=\vp_\pi:W_F\rightarrow\,^LM$ is the corresponding Langlands parameter for the $L$--packet $\Pi_M(\vp)$ containing $\pi.$  Let $R(\pi)$ be the Knapp-Stein $R$--group of $\pi$ and $R_{\vp,\pi}$ the Arthur $R$--group attached to $\vp$ and $\pi.$  Then $R(\pi)\simeq R_{\vp,\pi},$ and this isomorphism is realized by the map $\a\mapsto\check\a$ between roots and coroots.
\end{theorem}

\begin{proof}  By Lemma \ref{maxcase1} it is enough to prove this isomorphism in the case $\P$ is maximal.  This, however, is exactly the content of Corollary \ref{corollary siegel} and Theorem \ref{isomorphism max-non-siegel}. 
\end{proof}

\end{document}